\RequirePackage{fix-cm}
\documentclass{article}                     
\usepackage{graphicx}
\usepackage[utf8]{inputenc}
\usepackage{amsmath}
\usepackage{amssymb}
\usepackage [hmargin={3cm,3cm},vmargin={2.4cm,2.6cm}]{geometry}
\usepackage{hyperref}

\newtheorem{thm}{Theorem}[subsection]

\newtheorem{prop}{Proposition}[subsection]

\newtheorem{lem}{Lemma}[subsection]

\newtheorem{rem}{Remark}[subsection]

\newtheorem{cor}{Corollary}[subsection]

\newtheorem{defin}{Definition}[subsection]

\newtheorem{ex}{Example}[subsection]

\newtheorem{Ass}{Assumption}

\newtheorem{conj}{Conjecture}

\newcommand{\R}{\mathbb{R}}

\newcommand{\pr}{\textit{Proof : \vskip5mm}}

\newcommand{\cqfd}{\begin{flushright} 
                    $\Box$
                   \end{flushright}}

\usepackage{latexsym}
\begin{document}

\title{A Gaussian estimate for the heat kernel on differential forms and application to the Riesz transform}


\author{Baptiste Devyver         
}



\date{}

\maketitle

\begin{abstract}
Let $(M,g)$ be a complete Riemannian manifold which satisfies a Sobolev inequality of dimension $n$, and on which the volume growth is comparable to the one of $\R^n$ for big balls; if there is no non-zero $L^2$ harmonic 1-form, and the Ricci tensor is in $L^{\frac{n}{2}-\varepsilon}\cap L^\infty$ for an $\varepsilon>0$, then we prove a Gaussian estimate on the heat kernel of the Hodge Laplacian acting on 1-forms. This allows us to prove that, under the same hypotheses, the Riesz transform $d\Delta^{-1/2}$ is bounded on $L^p$ for all $1<p<\infty$. Then, in presence of non-zero $L^2$ harmonic 1-forms, we prove that the Riesz transform is still bounded on $L^p$ for all $1<p<n$, when $n>3$.
\end{abstract}

\section{Introduction and statements of the results}

\subsection{Riesz transform and heat kernel on differential forms}
Since Strichartz raised in 1983 the question whether the Riesz transform $d\Delta^{-1/2}$ is bounded on $L^p$ on a complete non-compact manifold $M$ (see \cite{Strichartz}), this problem has attracted a lot of attention. The litterature is too large to be cited extensively, but we refer the reader to the articles \cite{Coulhon-Duong1}, \cite{Coulhon-Duong2} and \cite{Auscher-Coulhon-Duong-Hofmann} for an overview of results in the field, as well as references. There have been several attempts to extend the classical Calder\'{o}n-Zygmund theory to the case of manifolds: the Calder\'{o}n-Zygmund decomposition argument which, in the case of $\R^n$, yields the $L^1\rightarrow L^1_w$ boundedness of the Riesz transform ($L^1_w$ being the weak $L^1$ space), has been adapted by Coulhon and Duong \cite{Coulhon-Duong2}. Let us begin with some definitions.

\begin{defin}

A Riemannian manifold $M$ has the \textbf{volume doubling property} if there exists a constant $C$ such that

\begin{equation}\tag{D}\label{doubling}
V(x,2R)\leq C\, V(x,R),\,\forall x\in M,\,\forall R>0,
\end{equation}
where $V(x,R)$ denotes the volume of the ball $B(x,R)$.

\end{defin}

\begin{defin}

Let us denote $p_t(x,y)$ the heat kernel of a complete Riemannian manifold $M$, i.e. the kernel of $e^{-t\Delta}$. We say that the heat kernel satisfies a \textbf{Gaussian upper-estimate} if there are some sonstants $c$ and $C$ such that

\begin{equation}\tag{G}\label{gaussian}
p_t(x,y)\leq \frac{C}{V(x,t^{1/2})}\exp\left(-c\frac{d^2(x,y)}{t}\right),\,\forall (x,y)\in M\times M,\,\forall t>0,
\end{equation}
where $d(x,y)$ denotes the distance between $x$ and $y$.

\end{defin}
With these definitions, the result of Coulhon and Duong states that if a complete Riemannian manifold has the volume doubling property \eqref{doubling} and its heat kernel satisfies a Gaussian upper-estimate \eqref{gaussian}, then the Riesz transform on $M$ is bounded on $L^p$ for every $1<p\leq 2$.\\
However, the duality argument which, in the case of $\R^n$, yields in turn that the Riesz transform is bounded on $L^p$ for $2\leq p<\infty$, does not apply in general for Riemannian manifolds. It has been observed that to get boundedness for the range $2\leq p<\infty$, it is enough to have a Gaussian estimate for the heat kernel of the Hodge Laplacian acting on differential 1-forms:

\begin{defin}

Let us denote $\vec{\Delta}=dd^\star+d^\star d$ the Hodge Laplacian acting on differential 1-forms, and $\vec{p_t}(x,y)$ the kernel of the corresponding heat operator $e^{-t\vec{\Delta}}$. If $x$ and $y$ are fixed, then $\vec{p_t}(x,y)$ is a linear morphism going from $\left(\Lambda^1T^\star M\right)_y$ to $\left(\Lambda^1T^\star M\right)_x$. We say that a \textbf{Gaussian estimate} holds for $\vec{p_t}(x,y)$ if there are some constants $C$, $c$ such that

\begin{equation}\tag{GF}\label{gaussian forms}
||\vec{p_t}(x,y)||\leq \frac{C}{V(x,t^{1/2})}\exp\left(-c\frac{d^2(x,y)}{t}\right),\,\forall (x,y)\in M\times M,\,\forall t>0,
\end{equation}
where $||.||$ is the operator norm.

\end{defin}
Then it is a consequence of the work of various authors \cite{Coulhon-Duong1}, \cite{Auscher-Coulhon-Duong-Hofmann}, \cite{Sikora} that if $M$ satisfies the volume doubling property \eqref{doubling} together with the Gaussian estimate for the heat kernel of the Hodge Laplacian on 1-forms \eqref{gaussian forms}, then the Riesz transform on $M$ is bounded on $L^p$ for all $2\leq p<\infty$. However, the only manifolds for which we know that \eqref{gaussian forms} holds are those with non-negative Ricci curvature (see \cite{Bakry} for example). 

Even outside the context of the Riesz transform problem, it is an interesting question to understand under which conditions \eqref{gaussian forms} holds on a complete non-compact manifold, since for example it implies the following pointwise gradient bound for the heat kernel on functions (see \cite{Coulhon-Duong1}):

$$|\nabla p_t(x,y)|\leq \frac{C}{\sqrt{t}V(x,t^{1/2})}\exp\left(-c\frac{d^2(x,y)}{t}\right),\,\forall (x,y)\in M\times M,\,\forall t>0,$$
which is one of the fundamental bounds (proved by Li-Yau  in \cite{Li-Yau}) for manifolds with non-negative Ricci curvature. In a recent paper, Coulhon and Zhang \cite{Coulhon-Zhang} have tried to prove bounds on the heat kernel of the Hodge Laplacian for a more general class of manifolds than those with non-negative Ricci curvature. To explain their result, we need some notations.\\
We denote by $Ric$ the symmetric endomorphism of $\Lambda^1T^\star M$ induced by the Ricci curvature tensor, and we can decompose it as

$$Ric=Ric_+-Ric_-,$$
where at a point $x$ of $M$, $(Ric_+)_x$ correspond to the positive eigenvalues of $Ric_x$, and $(Ric_-)_x$ to the non-positive eigenvalues. We then write $V(x)=\lambda^-(x)$, the lowest negative eigenvalue of $Ric_x$, with the convention that $V(x)=0$ if $Ric_x\geq 0$. The Bochner formula for the differential 1-forms writes

\begin{equation}\label{Bochner}\tag{B}
\vec{\Delta}=\bar{\Delta}+Ric=\bar{\Delta}+Ric_+-Ric_-,
\end{equation}
where $\bar{\Delta}=\nabla^\star\nabla$ is the so-called rough Laplacian. Then, if $p_t^V(x,y)$ denotes the kernel of $e^{-t(\Delta+V)}$, we have by the Kato inequality and domination theory (see \cite{H-S-U1}) 

$$||\vec{p_t}(x,y)||\leq |p_t^V(x,y)|,\,\forall (x,y)\in M\times M,\,\forall t>0.$$
Therefore Coulhon and Zhang, in order to bound $||\vec{p_t}(x,y)||$, study the heat kernel of the Schrödinger operator $\Delta+V$, making some assumptions on $V$. Contrary to $\vec{\Delta}$, $\Delta+V$ need not be non-negative, so they need to make some positivity assumption: Coulhon and Zhang make the assumption of \textbf{strong positivity} of $\Delta+V$: there is an $\varepsilon>0$ such that

$$\int_M V\varphi^2\leq (1-\varepsilon)\int_M |\nabla\varphi|^2,\,\forall \varphi\in C_0^\infty(M).$$
This is the translation at the quadratic forms level of the inequality

$$\Delta+V\geq \varepsilon \Delta.$$
The result of Coulhon and Zhang is the following:

\begin{thm}\label{Coulhon-Zhang}(Coulhon-Zhang, \cite{Coulhon-Zhang})\\
Let $M$ be a complete Riemannian manifold satisfying \eqref{doubling}, \eqref{gaussian}. We assume the non-collapsing of the volume of balls of radius $1$: for some constant $C>0$,

$$V(x,1)\geq C,\,\forall x\in M,$$
and that the negative part of the Ricci curvature is in $L^q\cap L^\infty$ for some $1\leq q<\infty$. We also assume that $\Delta+V$ is strongly positive. 

Then there is a constant $C$ such that

$$||\vec{p_t}(x,y)||\leq \min\left(\frac{t^\alpha}{V(x,t^{1/2})},1\right)\exp\left(-c\frac{d^2(x,y)}{t}\right),\,\forall (x,y)\in M\times M,\,\forall t>0,$$
where $\alpha$ is \textbf{strictly positive} and depends explicitely on $q$ and on $\varepsilon$; for example, if $q\geq 2$, one can take

$$\alpha=(q-1+\eta)(1-\varepsilon),$$
for all $\eta>0$.

\end{thm}
We want to make some comments about this result:

\begin{enumerate}

\item First, the estimate obtained on $\vec{p_t}(x,y)$ differs from the Gaussian estimate \eqref{gaussian forms} by a polynomial term in time $t^\alpha$. To explain where this extra polynomial term comes from, let us recall the method that Coulhon and Zhang develop in order to prove Theorem \ref{Coulhon-Zhang}. The main idea is to use a generalized \textit{Nash method}: the classical Nash method for the Laplacien shows that if $M$ satisfies a Sobolev inequality, then the heat kernel has on-diagonal estimates:

$$p_t(x,x)\leq \frac{C}{t^{n/2}},\,\forall t>0,\,\forall x\in M.$$
The classical Nash method uses in a crucial way the fact that $e^{-t\Delta}$ is contractive on $L^\infty$. Coulhon and Zhang generalize the Nash method to the case of a Schrödinger operator $\Delta+V$: the strong positivity assumption is the analogue of the Sobolev inequality of the classical Nash method, but this time the semi-group $e^{-t(\Delta+V)}$ is not in general contractive on $L^\infty$. Instead, Coulhon and Zhang show that under the hypotheses of Theorem \ref{Coulhon-Zhang},

$$||e^{-t(\Delta+V)}||_{\infty,\infty}\leq Ct^{\beta},\,\forall t>0,$$
for some $\beta>0$ (see Proposition 2.3 in \cite{Coulhon-Zhang}). Then, applying their generalized Nash method, they show an on-diagonal estimate for $p_t^V$, of the form $\frac{Ct^\alpha}{V(x,\sqrt{t})}$, where $\alpha$ is a function of $\beta$, and after they get the full off-diagonal estimates. \\
Let us emphasize that this extra polynomial term \textbf{does not} allow one to prove the boundedness of the Riesz transform on any $L^p$ space for $2\leq p<\infty$. 

\item Secondly, the geometric meaning of the strong positivity of $\Delta+V$ is not clear: to some extend, the strong positivity assumption is an hypothesis of triviality of the kernel of $\Delta+V$, yet the kernel of $\Delta+V$, contrary to the one of $\vec{\Delta}$, has no clear geometric meaning. It would be more natural to work directly with $\vec{\Delta}$, which has moreover the advantage of being non-negative.

\end{enumerate}
How to get rid of the extra polynomial term? As we mentionned in the first remark above, this extra polynomial term originates in the fact that instead of a uniform bound in time of $||e^{-t(\Delta+V)}||_{\infty,\infty}$, Coulhon and Zhang obtain an estimate of the form

$$||e^{-t(\Delta+V)}||_{\infty,\infty}\leq \frac{C}{t^\beta},\,\forall t>0,$$
where $\beta$ is some \textit{positive} constant. In the paper \cite{Simon}, B. Simon shows the following result for Schrödinger operators on $\R^n$, which improves on Coulhon and Zhang's estimate of $||e^{-t(\Delta+V)}||_{\infty,\infty}$ under slightly more restrictive hypotheses:

\begin{thm}\label{Simon}(B. Simon, \cite{Simon})\\
Let $V$ be a potential in $L^{\frac{n}{2}\pm\varepsilon}(\R^n)$ for some $\varepsilon>0$. We assume that $\Delta+V$ is strongly positive. Then 

$$||e^{-t(\Delta+V)}||_{\infty,\infty}\leq C,\,\forall t>0.$$

\end{thm} 
The proof consists in showing that there is a positive function $\eta$ in $L^\infty$, bounded from below by a positive constant, such that

$$(\Delta+V)\eta=0.$$
We want to show a similar result for the case of \textit{generalised Schrödinger operators}, that is operators of the form

$$L=\nabla^\star\nabla+\mathcal{R},$$
acting on the sections of a Riemannian vector bundle $E\rightarrow M$ endowed with a compatible connection $\nabla$, $\mathcal{R}$ being a field of symetric endomorphisms. A special case is the Hodge Laplacian acting on differential 1-forms. One of the main difficulties we will have to overcome in order to extend Theorem \ref{Simon} to this setting is that no notion of positivity is available for sections of a vector bundle. We will therefore not use the same approach as B. Simon.\\

\subsection{Our results} 
We will consider the class of manifolds satisfying a Sobolev inequality of dimension $n$:

\begin{defin}

We say that $M$ satisfies a Sobolev inequality of dimension $n$ if there is a constant $C$ such that 

\begin{equation}\tag{$\mathcal{S}_n$}\label{sobolev}
||f||_{\frac{2n}{n-2}}\leq C||\nabla f||_2,\,\forall f\in C_0^\infty(M).
\end{equation}

\end{defin}
\textit{Notation:} for two positive real function $f$ and $g$, we write

$$f\simeq g$$
if there are some positive constants $C_1$, $C_2$ such that 

$$C_1f\leq g\leq C_2 f.$$ 
In the following result (for the proof see Corollary \ref{estimee gaussienne sans cohomologie}), we extend the result of B. Simon to generalised Schrödinger operators in the class of manifolds satisfying a Sobolev inequality:

\begin{thm}\label{gaussian estimate}

Let $M$ be a complete, non-compact Riemannian manifold. We assume that $M$ satisfies a Sobolev inequality of dimension $n$ \eqref{sobolev}, and that the volume of big balls is Euclidean of dimension $n$:

$$V(x,R)\simeq R^n,\,\forall x\in M,\,\forall R\geq1.$$
Let $E\rightarrow M$ be a Riemannian vector bundle, endowed with a compatible connection $\nabla$, and let $L$ be a generalized Schrödinger operator acting on sections of $E$:

$$L=\nabla^\star\nabla+\mathcal{R},$$
$\mathcal{R}$ being a field of symmetric endomorphisms. Assume that $L$ is non-negative, that $\mathcal{R}_-$ lies in $L^{\frac{n}{2}-\varepsilon}\cap L^\infty$ for some $\varepsilon>0$, and that

$$Ker_{L^2}(L)=\{0\}.$$
Then the Gaussian estimate holds for $e^{-tL}$: if $K_{exp(-tL)}(x,y)$ denotes its kernel, for all $\delta>0$, there are two constants $C$ and $c$ such that

$$||K_{exp(-tL)}(x,y)||\leq \frac{C}{V(x,\sqrt{t})}\exp\left(-c\frac{d^2(x,y)}{t}\right),\,\forall (x,y)\in M\times M,\,\forall t>0.$$

\end{thm}
For the case where $L$ is the Hodge Laplacian acting on 1-forms, we get a boundedness result for the Riesz transform. Denote by $\mathcal{H}^1(M)$ the space of $L^2$ harmonic 1-forms:

$$\mathcal{H}^1(M)=\{\omega\in L^2(\Lambda^1 T^\star M) ; \vec{\Delta}\omega=0\}.$$
Then we show (cf Corollary \ref{transformee de Riesz bornee dans le cas sous-critique}):

\begin{thm}\label{Riesz 1<p<infty}

Let $M$ be a complete, non-compact Riemannian manifold, of dimension $m$. We assume that $M$ satisfies the Sobolev inequality of dimension $n$ \eqref{sobolev}, and that the volume of big balls is Euclidean of dimension $n$:

$$V(x,R)\simeq R^n,\,\forall x\in M,\,\forall R\geq1.$$
Assume that the negative part of the Ricci curvature is in $L^{\frac{n}{2}-\varepsilon}\cap L^\infty$ for some $\varepsilon>0$, and that

$$\mathcal{H}^1(M)=\{0\}.$$
Then the Gaussian estimate holds for $e^{-t\vec{\Delta}}$: for all $\delta>0$, there are two constants $C$ and $c$ such that

$$||\vec{p_t}(x,y)||\leq \frac{C}{V(x,\sqrt{t})}\exp\left(-c\frac{d^2(x,y)}{t}\right),\,\forall (x,y)\in M\times M,\,\forall t>0,$$
$\vec{p_t}(x,y)$ denoting the kernel of $e^{-t\vec{\Delta}}$. Moreover, the Riesz transform on $M$ is bounded on $L^p$ for all $1<p<\infty$.

\end{thm}

\begin{rem}

The hypothesis $\mathcal{H}^1(M)=\{0\}$ is somewhat optimal to get the boundedness of the Riesz transform on $L^p$ for $1<p<\infty$, in the class of manifolds considered: it is known that the Riesz transform on the connected sum of two euclidean spaces $\R^n\#\R^n$ is bounded on $L^p$ if and only if $p\in (1,n)$ (cf \cite{Carron-Coulhon-Hassell}). And $\R^n\#\R^n$ satisfies all the hypotheses of Theorem \ref{Riesz 1<p<infty}, except that $\mathcal{H}^1(\R^n\#\R^n)\neq\{0\}$: indeed, if $n\geq3$, due to the fact that $\R^n\#\R^n$ has two non-parabolic ends, we can find a non-constant harmonic function $h$ such that $\nabla h$ is $L^2$. Then $dh$ is a non-zero, $L^2$ harmonic 1-form. We conjecture that in fact, the existence of several ends is the only obstruction to the boundedness of the Riesz transform in the range $(1,\infty)$ in the class of manifolds that we consider:

\end{rem}

\begin{conj}

Let $M$ be a complete, non-compact Riemannian manifold, of dimension $m$. We assume that $M$ satisfies a Sobolev inequality of dimension $n$ \eqref{sobolev}, and that the volume of big balls is Euclidean of dimension $n$:

$$V(x,R)\simeq R^n,\,\forall x\in M,\,\forall R\geq1.$$
Assume that the negative part of the Ricci curvature is in $L^{\frac{n}{2}-\varepsilon}\cap L^\infty$ for some $\varepsilon>0$, and that $M$ has only one end. Then the Riesz transform on $M$ is bounded on $L^p$ for every $1<p<\infty$.

\end{conj}
In the second part of the paper, we adress the following question: what happens for the Riesz transform if one removes the hypothesis that $\mathcal{H}^1(M)=\{0\}$ ? As we see from the example of $\R^n\#\R^n$, all we can hope is the boundedness on $L^p$ for $1<p<n$. And indeed, we will show (cf Theorem \ref{transformee de Riesz pour 1<p<n}):

\begin{thm}\label{Riesz 1,n}

Let $M$ be a complete, non-compact Riemannian manifold. For an $n>3$, we assume that $M$ satisfies a Sobolev inequality of dimension $n$ \eqref{sobolev}, and that the volume of big balls is Euclidean of dimension $n$:

$$V(x,R)\simeq R^n,\,\forall x\in M,\,\forall R\geq1.$$
Assume also that the negative part of the Ricci curvature is in $L^{\frac{n}{2}-\varepsilon}\cap L^\infty$ for some $\varepsilon>0$. Then the Riesz transform on $M$ is bounded on $L^p$, for all $1<p<n$.

\end{thm}
The strategy of the proof is a perturbation argument, using the Gaussian estimate proved in Theorem \ref{gaussian estimate}: first, we know from our hypotheses and \cite{Coulhon-Duong2} that the Riesz transform on $M$ is bounded on $L^p$ for all $1<p\leq 2$. From the proof of Theorem \ref{Riesz 1<p<infty}, there is a constant $\eta$ such that if 

$$||Ric_-||_{L^{n/2}}\leq \eta,$$
then $Ker_{L^2}(\vec{\Delta})=\{0\}$. We take $V$ a smooth, non-negative, compactly supported potential such that 

$$||(Ric+V)_-||_{L^{n/2}}\leq \eta,$$
then we will have (by Theorem \ref{gaussian estimate}) a Gaussian estimate for $e^{-t(\vec{\Delta}+V)}$. We will see that this implies the boundedness of the Riesz transform with potential $d(\Delta+V)^{-1/2}$ on $L^p$, for all $\frac{n}{n-1}<p<n$. Then, using a perturbation argument from \cite{Carron4}, we will be able to show that $d\Delta^{-1/2}-d(\Delta+V)^{-1/2}$ is bounded on $L^p$ for all $\frac{n}{n-1}<p<n$.\\\\
We will also study the $L^p$ reduced cohomology: we let $H^1_p(M)$, the first space of $L^p$ \textbf{reduced cohomology}, to be the quotient of $\{\alpha\in L^p : d\alpha=0\}$ by the closure in $L^p$ of $dC_0^\infty(M)$. We have the following result which sums up the results of this paper, including the one concerning the $L^p$ reduced cohomology:

\begin{thm}

Let $M$ be a complete, non-compact Riemannian manifold. We assume that $M$ satisfies a Sobolev inequality of dimension $n$ \eqref{sobolev}, and that the volume of big balls is Euclidean of dimension $n$:

$$V(x,R)\simeq R^n,\,\forall x\in M,\,\forall R\geq1.$$
Assume that the negative part of the Ricci curvature is in $L^{\frac{n}{2}-\varepsilon}\cap L^\infty$ for some $\varepsilon>0$. Then we have the following alternative:

\begin{enumerate}
	\item $\mathcal{H}^1(M)$, the space of $L^2$ harmonic 1-forms, is trivial. Then for all $1<p<\infty$, the Riesz transform on $M$ is bounded on $L^p$, and $H^1_p(M)$, the first space of $L^p$ reduced cohomology of $M$, is trivial.
	\item $\mathcal{H}^1(M)$ is not trivial. If $n>3$, then for all $1<p<n$, the Riesz transform on $M$ is bounded on $L^p$ and $H_p^1(M)\simeq \mathcal{H}^1(M)$. Moreover, if $M$ has more than one end, for $p\geq n$ the Riesz transform on $M$ is not bounded on $L^p$, and $H_p^{1}(M)$ is not isomorphic to $\mathcal{H}^1(M)$.
\end{enumerate}

\end{thm}

\section{Preliminaries}

We consider an operator $L$ of the form $\nabla^*\nabla + \mathcal{R}_+-\mathcal{R}_-$, acting on a Riemannian fiber bundle $E\rightarrow M$, where $\nabla$ is a connection on $E\rightarrow M$ compatible with the metric, and for $p\in M$, $\mathcal{R}_+(p)$, $\mathcal{R}_-(p)$ are non-negative symmetric endomorphism acting on the fiber $E_p$. We will say that $L$ is a \textbf{generalised Schr\"{o}dinger operator}. Let us denote $\bar{\Delta}:=\nabla^*\nabla$, the "rough Laplacian", and $C^\infty(E)$ (resp. $C_0^\infty(E)$) the set of smooth sections of $E$ (resp., of smooth sections of $E$ which coincide with the zero section outside a compact set).\\
We define $H:=\bar{\Delta}+\mathcal{R}_+$. We will consider the $L^2$-norm on sections of $E$:
$$||\omega||_2^2=\int_M|\omega|^2(p)dvol(p),$$
where $|\omega|(p)$ is the norm of the evaluation of $\omega$ in $p$ and $dvol$ the Riemannian volume form. We will denote $L^2(E)$, or simply $L^2$ when there is no confusion possible, for the set of sections of $E$ with finite $L^2$-norm.\\
We have in mind the case of $\vec{\Delta}=d^*d+dd^*$, the Hodge Laplacian, acting on differential forms of degree $1$, for which we have the \textbf{Bochner decomposition}:
$$\vec{\Delta}=\bar{\Delta}+Ric,$$
where $\bar{\Delta}=\nabla^*\nabla$ is the rough Laplacian on 1-forms, and $Ric\in End(\Lambda^1T^*M)$ is canonically identified -- using the metric -- with the Ricci tensor.\\
From classical results in spectral analysis (an obvious adaptation to $\bar{\Delta}$ of Strichartz's proof that the Laplacian is self-adjoint on a complete manifold, see Theorem 3.13 in \cite{Pigola-Rigoli-Setti}), we know that if $\mathcal{R}_-$ is bounded, then $\bar{\Delta}+\mathcal{R}_+-\mathcal{R}_-$ is essentially self-adjoint on $C_0^\infty(E)$.

\subsection{Consequences of the Sobolev inequality}

Let us denote $q_H$ the quadratic form associated to $H$:

$$q_H(\omega)=\int_M|\nabla\omega|^2+\int_M\langle \mathcal{R}_+\omega,\omega\rangle,$$
and $q_\Delta$ the quadratic form associated to the usual Laplacian on functions:

$$q_\Delta(u)=\int_M|du|^2.$$
We will see in this section that $H=\bar{\Delta}+\mathcal{R}_+$ shares with the usual Laplacian acting on functions a certain amount of functionnal properties. It is due to the following \textit{domination property}, consequence of the Kato inequality (see \cite{Berard-Besson}):

\begin{prop}\label{domination}[Domination property]
 
For all $\omega\in C_0^\infty(E)$,
$$|e^{-tH}\omega|\leq e^{-t\Delta}|\omega|,\,\forall t\geq0,$$
and

$$|H^{-\alpha}\omega|\leq \Delta^{-\alpha}|\omega|,\,\forall\alpha>0.$$

\end{prop}
\pr The first part comes directly from \cite{Berard-Besson}. The second domination is a consequence of the first one and of the following formulae:
$$H^{-\alpha}=\frac{1}{\Gamma(\alpha+1)}\int_0^\infty e^{-tH}t^{\alpha-1}dt,$$
$$\Delta^{-\alpha}=\frac{1}{\Gamma(\alpha+1)}\int_0^\infty e^{-t\Delta}t^{\alpha-1}dt.$$
\cqfd
The Kato inequality also yields:

\begin{prop}\label{Sobolev pour H}
$H$ satisfies a Sobolev inequality of dimension $n$: there is a constant $C$ such that
$$||\omega||^2_{\frac{2n}{n-2}}\leq C\langle H\omega,\omega\rangle,\,\forall\omega\in C_0^\infty(E).$$
\end{prop}
The domination property (Proposition \ref{domination}), together with the fact that $e^{-t\Delta}$ is a contraction semigroup on $L^p$ for all $1\leq p\leq\infty$, gives at once that $e^{-tH}$ is also a contraction semigroup on all the $L^p$ spaces. By Stein's Theorem (Theorem 1 p.67 in \cite{Stein}), $e^{-tH}$ is analytic bounded on $L^p$, for all $1<p<\infty$. \\
From the ultracontractivity estimate: 

$$||e^{-t\Delta}||_{1,\infty}\leq \frac{C}{t^{n/2}},\,\forall t>0,$$
valid since $M$ satisfies a Sobolev inequality of dimension $n$ (see \cite{Saloff-Coste1}), and the domination of Proposition \ref{domination}, we deduce that we also have:
$$||e^{-tH}||_{1,\infty}\leq \frac{C}{t^{n/2}},\,\forall t>0.$$
By interpolation with $||e^{-tH}||_{\infty,\infty}\leq 1$, we deduce that for all $1\leq p\leq \infty$, there exists $C$ such that 
$$||e^{-tH}||_{p,\infty}\leq\frac{C}{t^{n/2p}},\,\forall t>0.$$
Interpolating with $||e^{-tH}||_{p,p}\leq 1$, we obtain:

$$||e^{-tH}||_{p,q}\leq\frac{C}{t^{\frac{n}{2p}\left( 1-\frac{p}{q}\right)}}.$$
To sum up, we have proved:

\begin{cor}\label{e^-tH semi-groupe de contraction}
$e^{-tH}$ is a contraction semigroup on $L^p$, for all $1\leq p\leq \infty$.\\
For all $1\leq p\leq \infty$, there exists $C$ such that:

$$||e^{-tH}||_{p,q}\leq\frac{C}{t^{\frac{n}{2p}\left( 1-\frac{p}{q}\right)}},\,\forall t>0,\,\forall q>p.$$
Moreover, $e^{-tH}$ is analytic bounded on $L^p$ with sector of angle $\frac{\pi}{2}\left(1-\left|\frac{2}{p}-1\right|\right)$, for all $1<p<\infty$.
\end{cor}
We recall the following consequences of the analyticity of a semigroup, which come from the Dunford-Schwarz functionnal calculus (see \cite{Reed-Simon2}, p.249):

\begin{cor}\label{analyticity}

Let $e^{-zA}$ an analytic semigroup on a Banach space $X$. Then there exists a constant $C$ such that for all $\alpha>0$:

\begin{enumerate}
 \item $$||A^\alpha e^{-tA}||\leq \frac{C}{t^\alpha},\,\forall t>0.$$
 \item $$||(I+tA)^\alpha e^{-tA}||\leq C,\,\forall t>0.$$
\end{enumerate}

\end{cor}
Furthermore, we also have the following simple but important consequences of the domination:

\begin{thm}\label{inegalites fonctionnelles}
 
$H$ satisfies the following properties:

\begin{enumerate}
 \item The \textbf{mapping properties}:\\ 
 For all $\alpha>0$, 
 
 $$H^{-\alpha/2} : L^p\longrightarrow L^q$$ 
 is bounded whenever $\frac{1}{q}=\frac{1}{p}-\frac{\alpha}{n}$ and $p<q<\infty$ (in particular we must have $p<\frac{n}{\alpha}$).
 
 \item The \textbf{Gagliardo-Nirenberg inequalities}:\\ 
 For all $s\geq r>n$,

$$||\omega||_{L^\infty}\leq C(n,r,s)||H\omega||^\theta_{r/2}||\omega||^{1-\theta}_{s/2},\,\forall \omega\in C_0^\infty(E),$$
where $\theta=\frac{n/s}{1-(n/r)+(n/s)}$.
 
\end{enumerate}

\end{thm}
\pr The mapping properties for $H$ are the consequence of the domination of Proposition \ref{domination} and of the mapping properties for $\Delta$, which hold since $M$ satisfies a Sobolev inequality of dimension $n$ (cf \cite{Varopoulos}, Theorem 1 and \cite{Coulhon-Saloff-Varopoulos}, Theorem II.4.1). The Gagliardo-Nirenberg inequalities are extracted from \cite{Coulhon}, Theorems 1 and 2, given the ultracontractivity of $e^{-tH}$ and its analiticity on $L^p$ for $1<p<\infty$ (Corollary \ref{e^-tH semi-groupe de contraction}).
\cqfd
Finally, we have the following important fact:

\begin{prop}\label{independance of lambda}

All the results of this section are also valid if we replace $H$ by $H+\lambda$ with $\lambda>0$, and moreover the constants in the Sobolev inequality, in the Gagliardo-Nirenberg inequality and also the norms of the operators $(H+\lambda)^{-\alpha} : L^p\rightarrow L^q$, are all bounded \textbf{independantly} of $\lambda$. 

\end{prop}
This will be used intensively later.

\subsection{Strong positivity}

As in the previous section, denote $H:=\bar{\Delta}+\mathcal{R}_+$ and 

$$L=\nabla^\star\nabla+\mathcal{R}=H-\mathcal{R}_-.$$

We assume -- as it is the case for the Hodge Laplacian acting on 1-forms -- that $L$ is a \textit{non-negative} operator:
 
\begin{Ass}\label{operateur positif}
 
$L$ is a \textit{non-negative} operator.

\end{Ass}
It is equivalent to the following inequality: if $\omega\in C_0^\infty(E)$, 

$$0\leq\langle \mathcal{R}_-\omega,\omega\rangle\leq \langle H\omega,\omega\rangle.$$
Let us recall the following classical definition:
\begin{defin}
 
The Hilbert space $H_0^1$ is the completion of $C_0^\infty(E)$ for the norm given by the quadratic form associated to the self-adjoint operator $H$. 
\end{defin}
We also recall some of the properties of this space $H_0^1$ associated to $H$:

\begin{prop}\label{non-parabolique}
\begin{enumerate}

 \item \label{enu1} $H_0^1\hookrightarrow L^{\frac{2n}{n-2}}(E)$. In particular, it is a space of sections of $E\rightarrow M$.
 
 \item \label{enu2} $H^{1/2}$, defined on $C_0^\infty(E)$, extends uniquely to a bijective isometry from $H_0^1$ to $L^2(E)$. \\ 
 Thus we can consider $H^{-1/2} : L^2(E)\rightarrow H_0^1$.
 
 \item \label{enu3} If we consider the operator $H^{1/2}$ given by the Spectral Theorem - denote it $H_s^{1/2}$ to avoid confusion with the one we have just defined from $H_0^1$ to $L^2$ - then $\mathcal{D}om(H_s^{1/2})=H_0^1\cap L^2(E)$, and moreover $H^{1/2}$ coincide with $H_s^{1/2}$ on $H_0^1\cap L^2(E)$. 
 
\end{enumerate}

\end{prop}
\textit{Sketch of proof of Proposition (\ref{non-parabolique}):}\\\\
(\ref{enu1}) is a consequence of the Sobolev inequality of Proposition \ref{Sobolev pour H}. The Sobolev inequality implies that $H$ is non-parabolic, and (\ref{enu2}) can be obtained by the same method as in \cite{Devyver}. (\ref{enu3}) can also be obtained by the techniques developped in \cite{Devyver} in the context of Schr\"{o}dinger operators acting on functions, which adapts to the case of Schr\"{o}dinger operators acting on sections of a vector bundle.
\cqfd
In what follows, we assume that $\mathcal{R}_-\in L^{\frac{n}{2}}$.

\begin{defin}\label{sous-criticite}
 
We say that $L$ is \textbf{strongly positive} if one of the following equivalent -- at least when $\mathcal{R}_-\in L^{\frac{n}{2}}$ -- conditions is satisfied :

\begin{enumerate}

 \item There exists $\varepsilon>0$ such that:
 
 $$0\leq\langle \mathcal{R}_-\omega,\omega\rangle\leq (1-\varepsilon) \langle H\omega,\omega\rangle,\,\forall\omega\in C_0^\infty(E).$$
 
 \item $$Ker_{H_0^1}(L)=\{0\}.$$ 
 
 \item The (non-negative, self-adjoint compact if $\mathcal{R}_-\in L^{\frac{n}{2}}$) operator $A:=H^{-1/2}\mathcal{R}_-H^{-1/2}$ acting on $L^2(E)$ satisfies:
 
 $$||A||_{2,2}\leq (1-\varepsilon),$$ 
 where $\varepsilon>0$. 
 
\end{enumerate}

\begin{rem}
 
In general, we have the equivalence between 1) and 3) and the implication 3)$\Rightarrow$2), under the sole hypothesis that $A$ is self-adjoint (which is the case if $\mathcal{R}_-\in L^{\frac{n}{2}}$, but can be true under more general conditions). The implication 2)$\Rightarrow$3) is true as soon as $A$ is self-adjoint compact.

\end{rem}

\end{defin}
\textit{Proof of the equivalence:}\\\\
We can write:

$$L=H-\mathcal{R}_-=H^{1/2}(I-A)H^{1/2}.$$
First, let us prove that $1)\Leftrightarrow 3')$, where $3')$ is defined to be:

$$3'):\,\langle Au,u\rangle \leq (1-\varepsilon)\langle u,u\rangle,\,\forall u\in L^2.$$
Remark that $3')$ is equivalent to $3)$ when $A$ is self-adjoint. Let $\omega\in C_0^\infty(E)$, and set $u=H^{1/2}\omega\in L^2(E)$. Then 

$$\begin{array}{rcl}
 
\langle Au,u\rangle \leq (1-\varepsilon)\langle u,u\rangle &\Leftrightarrow& \langle H^{-1/2}\mathcal{R}_-\omega,H^{1/2}\omega\rangle\leq (1-\varepsilon) \langle H^{1/2}\omega,H^{1/2}\omega\rangle\\

&\Leftrightarrow& \langle \mathcal{R}_-\omega,\omega\rangle\leq (1-\varepsilon)\langle H\omega,\omega\rangle

\end{array}$$
Then we show that $3)\Rightarrow 2)$. This is a consequence of the following Lemma:

\begin{lem}\label{isometrie entre noyaux}

If $A$ is self-adjoint, then
 
$$H^{1/2} : Ker_{H_0^1}(L)\rightarrow Ker_{L^2}(I-A)$$
is an isomorphism (and it is of course an isometry).

\end{lem}
\pr Let $u\in H_0^1$; we can write  $u=H^{-1/2}\varphi$, where $\varphi\in L^2(E)$. By definition, $Lu=0$ means that for every $v\in C_0^\infty(E)$,\\

$$\langle u,Lv\rangle=0.$$
This equality makes sense, because since $H$ satisfies a Sobolev inequality, $H_0^1\hookrightarrow L_{loc}^1$. The Spectral Theorem then implies, since $C_0^\infty\subset\mathcal{D}om(H)$, that given $v\in C_0^\infty(E)$ the following equality holds in $L^2(E)$:

$$Hv=H^{1/2}H^{1/2}v.$$
Hence

$$Lv=(H-R_-)v=H^{1/2}(I-A)H^{1/2}v.$$
Let $w:=(I-A)H^{1/2}v$; then the preceeding equality shows that $w\in \mathcal{D}om(H^{1/2})=H_0^1\cap L^2(E)$. Furthermore, $H^{1/2}w=Hv$ is compactly supported, so we have:

$$\langle u,H^{1/2}w\rangle=\langle H^{1/2}u,w\rangle.$$
Indeed, if $u\in H_0^1\cap L^2$ it is a consequence of Lemma 3.1 in \cite{Devyver}, and a limiting argument plus the fact that $H_0^1\hookrightarrow L^2_{loc}$ shows that it is true for all $u\in H_0^1$.

$$Lu=0\Longleftrightarrow \forall v\in C_0^\infty,\,\langle H^{1/2}u,(I-A)H^{1/2}v\rangle=0.$$
But since $H^{1/2}C_0^\infty(E)$ is dense in $L^2(E)$, we get, using the fact that $A$ is self-adjoint:

$$\begin{array}{rcl}
Lu=0&\Longleftrightarrow &\forall v\in L^2,\,\langle H^{1/2}u,(I-A)v\rangle=0\\\\
&\Longleftrightarrow& H^{1/2}u\in Ker_{L^2}(I-A)
\end{array}$$
\cqfd
It remains to prove that $2)\Rightarrow 3)$ ; this is a consequence of Lemma \ref{isometrie entre noyaux} and of the following Lemma, which is extracted from Proposition 1.2 in \cite{Carron1}:

\begin{lem}\label{cas Ric dans L^n/2}
 
Assume $\mathcal{R}_-\in L^{\frac{n}{2}}$. Then $A:=H^{-1/2}\mathcal{R}_-H^{-1/2}$ is a non-negative, self-adjoint compact operator on $L^2(E)$. Moreover, 

$$||A||_{2,2}\leq C||\mathcal{R}_-||_{\frac{n}{2}},$$
where $C$ depends only on a (any) Sobolev constant $K$ for $H$, that is any constant $K$ such that

$$||\omega||^2_{\frac{2n}{n-2}}\leq K \langle H\omega,\omega\rangle.$$

\end{lem}
\cqfd
We will also need the following Lemma, which is an easy consequence of the definition of strong positivity:

\begin{lem}\label{Sobolev pour sous-critique}
 
Let $H$ be of the form: $H=\bar{\Delta}+\mathcal{R}_+$, with $\mathcal{R}_+$ non-negative. Let $\mathcal{R}_-\in End(\Lambda^1 T^*M)$ be symmetric, non-negative, in $L^{\frac{n}{2}}$ such that $L:=H-\mathcal{R}_-$ is strongly positive. Then a Sobolev inequality of dimension $n$ is valid for $L$, i.e.

$$||\omega||^2_{\frac{2n}{n-2}}\leq C\langle L\omega,\omega\rangle, \,\forall\omega\in C_0^\infty(E).$$ 

\end{lem}
\pr By definition of strong positivity,

$$\langle \mathcal{R}_-\omega,\omega\rangle\leq (1-\varepsilon)\langle H\omega,\omega\rangle.$$
Therefore:

$$\begin{array}{rcl}
 
\langle L\omega,\omega\rangle&=&\langle H\omega,\omega\rangle-\langle \mathcal{R}_-\omega,\omega\rangle \\\\

&\geq&(1-(1-\varepsilon))\langle H\omega,\omega\rangle \\\\

&\geq&\varepsilon C||\omega||^2_{\frac{2n}{n-2}},

\end{array}$$
where we have used in the last inequality the fact that $H$ satisfies a Sobolev inequality.
\cqfd

\section{Gaussian upper-bound for the Heat Kernel on 1-forms}

\subsection{Estimates on the resolvent of the Schr\"{o}dinger-type operator}

In this section, we will show how to obtain bounds on the resolvent of $L:=\nabla^*\nabla+\mathcal{R}_+-\mathcal{R}_-=H-\mathcal{R}_-$. In order to do this, we first have to estimate the resolvent of the operator $H=\bar{\Delta}+\mathcal{R}_+$. Recall from Corollary \ref{e^-tH semi-groupe de contraction} that $e^{-tH}$ is a contraction semigroup on $L^p$, for all $1\leq p\leq \infty$. Using the formula:
$$(L+\lambda)^{-1}=\int_0^\infty e^{-tL}e^{-t\lambda}dt,$$
we get:

\begin{prop}\label{resolvante de H}
 
For all $\lambda>0$ and for all $1\leq p\leq \infty$,

$$||(H+\lambda)^{-1}||_{p,p}\leq  \frac{1}{\lambda}.$$

\end{prop}

\begin{rem}
 
The case $p=\infty$ is by duality, for $(H+\lambda)^{-1}$ is defined on $L^\infty$ by duality. Indeed, for $g\in L^\infty$, we define $(H+\lambda)^{-1}g$ so that:

$$\langle (H+\lambda)^{-1}g,f\rangle:=\langle g,(H+\lambda)^{-1}f\rangle,\,\forall f\in L^1$$
(recall that $(L^1)^\prime=L^\infty$). It is then easy to see that $||(H+\lambda)^{-1}||_{1,1}\leq  \frac{1}{\lambda}$ implies $||(H+\lambda)^{-1}||_{\infty,\infty}\leq  \frac{1}{\lambda}$.

\end{rem}
We now estimate the resolvent of the operator $L:=\nabla^*\nabla+\mathcal{R}_+-\mathcal{R}_-$; as before, $L$ acts on the sections of a vector bundle $E\rightarrow M$ (see the beginning of the Preliminaries for the general context). The key result is the following:

\begin{thm}\label{resolvante du Laplacien de Hodge}
 
Let $(M,g)$ be a complete Riemannian manifold which satisfies a Sobolev inequality of dimension $n$ \eqref{sobolev}, and suppose that $\mathcal{R}_-$ is in $L^{\frac{n}{2}\pm\varepsilon}$ for some $\varepsilon>0$. We also assume that $L$, acting on the sections of $E\rightarrow M$, is strongly positive. Then for all $1\leq p\leq\infty$, there exists a constant $C(p)$ such that 

$$||(L+\lambda)^{-1}||_{p,p}\leq \frac{C(p)}{\lambda},\,\forall \lambda>0.$$

\end{thm}
\pr In this proof, we write $L^q$ for $L^q(E)$. Let us denote $H_\lambda:=H+\lambda$. So

$$(L+\lambda)^{-1}=(I-T_\lambda)^{-1}H_\lambda^{-1},$$
where $T_\lambda:=H_\lambda^{-1}R_-$. We will prove that $(I-T_\lambda)^{-1}$ is a bounded operator on $L^p$, with norm independant of $\lambda$, then applying Proposition \ref{resolvante de H} gives the result. To achieve this, we will show that the series $\sum_{n\geq 0} T_\lambda^n$ converges in $\mathcal{L}(L^p,L^p)$, uniformly with respect to $\lambda\geq0$.\\
The aim of the next two Lemmas is to prove that $T_\lambda$ acts on all the $L^q$ spaces. We single out the case $q=\infty$, for it requires a different ingredient for its proof:

\begin{lem}\label{estimee infinie}
 
$T_\lambda : L^\infty\longrightarrow L^\infty$ is bounded as a linear operator, uniformly with respect to $\lambda\geq0$.

\end{lem}
\pr We have seen that $e^{-tH_\lambda}$ satisfies the mapping properties and the Gagliardo-Nirenberg inequalities of Theorem \ref{inegalites fonctionnelles} with constants independant of $\lambda\geq0$. Let $u\in L^\infty$. We apply the Gagliardo-Nirenberg inequality for $H_\lambda$: 

$$||T_\lambda u||_\infty\leq C||R_- u||^\theta_{n/2+\varepsilon}||T_\lambda u||_p^{1-\theta},\,\forall p,$$
with $C$ independant of $\lambda$ (see Proposition \ref{independance of lambda}). We have $||R_-u||_{\frac{n}{2}+\varepsilon}\leq ||R_-||_{\frac{n}{2}+\varepsilon }||u||_\infty$. By the mapping properties of Theorem \ref{inegalites fonctionnelles}, $H_\lambda^{-1} : L^{\frac{n}{2}-\varepsilon}\rightarrow L^s$ for a certain $s$, with a norm bounded independantly of $\lambda$. So we get:

$$||T_\lambda u||_\infty\leq C||R_-||^\theta_{n/2+\epsilon}(||H_\lambda^{-1}||_{n/2-\varepsilon,s}||R_-||_{n/2-\varepsilon})^{1-\theta}||u||_\infty\leq C||u||_\infty$$
\cqfd

\begin{lem}\label{action de T sur les L^p}

\begin{enumerate}

\item For all $1\leq \beta\leq\infty$,

$$\mathcal{R}_- : L^{\beta}\rightarrow L^{\frac{n\beta}{n+2\beta}}$$
is bounded.
 
\item There exists $\nu>0$ (small and independant of $\lambda\geq0$), such that for all $\beta<\infty$, and for all $\lambda\geq 0$,

$$T_\lambda : L^\beta\rightarrow L^{r}\cap L^{s},$$
where $\frac{1}{r}=\max(\frac{1}{\beta}-\nu,0^+))$ and $\frac{1}{s}=\min(\frac{1}{\beta}+\nu,1)$, is bounded uniformly with respect to $\lambda$ (here $0^+$ denotes any positive number).

\item For $\beta=\infty$,

$$T_\lambda : L^\infty\rightarrow L^\infty\cap L^p$$
is bounded uniformly with respect to $\lambda$, if $p$ big enough.

\item For $\beta$ large enough,

$$T_\lambda : L^\beta\rightarrow L^\beta\cap L^\infty$$
is bounded uniformly with respect to $\lambda$.

\end{enumerate}

\end{lem}
\pr If $u\in L^\beta$ and $v\in L^\gamma$, then

$$||uv||_{\frac{\gamma\beta}{\gamma+\beta}}\leq ||u||_\beta||v||_\gamma.$$
Therefore, $\mathcal{R}_- : L^\beta\rightarrow L^q$ is bounded, where $\frac{1}{q}=\frac{1}{\beta}+\frac{1}{p}$, for all $p\in[\frac{n}{2}-\varepsilon,\frac{n}{2}+\varepsilon]$. Taking $p=\frac{n}{2}$, we find the first result of the Lemma.\\
Applying the mapping property of Theorem \ref{inegalites fonctionnelles}, we deduce that:

$$H_\lambda^{-1}\mathcal{R}_- : L^\beta\rightarrow \rightarrow L^r\cap L^s$$
is bounded independantly of $\beta$, and also uniformly with respect to $\lambda\geq0$ by Proposition \ref{independance of lambda}, where 

$$\frac{1}{r}=\max\left(\left(\frac{2}{n+2\varepsilon}-\frac{2}{n}\right)+\frac{1}{\beta},0^+\right)=\max\left(\frac{1}{\beta}-\mu,0^+\right),$$
and 
$$\frac{1}{s}=\min\left(\left(\frac{2}{n-2\varepsilon}-\frac{2}{n}\right)+\frac{1}{\beta},1\right)=\min\left(\frac{1}{\beta}+\mu',1\right),$$
hence the second part of the Lemma with $\nu=\min(\mu,\mu')$.\\\\
For the case $\beta=\infty$, we have $s=\frac{1}{\mu'}=p$ large, and we already know from Lemma \ref{estimee infinie} that $T_\lambda$ send $L^\infty$ to $L^\infty$.\\\\
For the case $\beta$ large enough: since $\mathcal{R}_-\in L^{\frac{n}{2}+\varepsilon}$, if $\beta$ is large enough and $u\in L^\beta$, then $\mathcal{R}_-u\in L^{\frac{n}{2}+\alpha}$ for an $\alpha>0$. We apply the Gagliardo-Nirenberg inequality: for such a $\beta$,

$$||H_\lambda^{-1}\mathcal{R}_-u||_\infty\leq C||\mathcal{R}_-u||^\theta_{n/2+\alpha}||H_\lambda^{-1}\mathcal{R}_-u||^{1-\theta}_\beta.$$
This yields the result.
\cqfd
As a corollary of Lemma \ref{action de T sur les L^p}, we obtain:

\begin{prop}\label{action de T^N sur les L^p}
 
For all $1\leq\beta\leq \infty$ and $1\leq \alpha\leq\infty$, there exists an $N\in \mathbb{N}$ (which depends only on $\beta$ and $\alpha$, and not on $\lambda$), such that for all $\lambda\geq0$,

$$T_\lambda^N : L^\alpha\rightarrow L^\beta$$
is bounded uniformly with respect to $\lambda$.

\end{prop}
Now, we have the following

\begin{lem}\label{action de T sur L^beta}
 
Let $\beta:=\frac{2n}{n-2}$. Then $||T_\lambda^k||_{\beta,\beta}\leq C(1-\mu)^k$ for all $k\in\mathbb{N}$ with constants $C$ and $0<\mu<1$ independant of $\lambda\geq0$.
\end{lem}
Let us postpone the proof of Lemma \ref{action de T sur L^beta}. Let us show that it implies the convergence of the series $\sum_{n\geq 0} T_\lambda^n$ in $\mathcal{L}(L^p,L^p)$ for all $1\leq p\leq\infty$, uniformly with respect to $\lambda\geq0$. Indeed, for a fixed $p$, according to Proposition \ref{action de T^N sur les L^p} we can find $N$ such that

$$T_\lambda : L^\beta\rightarrow L^p$$
and

$$T_\lambda : L^p\rightarrow L^\beta.$$
Then, for $k\geq 2N$, we decompose $T_\lambda^k$ in

$$T_\lambda^k=T_\lambda^NT_\lambda^{N-2k}T_\lambda^N,$$
where $T_\lambda^N$ on the right goes from $L^p$ to $L^\beta$, and $T_\lambda^N$ on the left goes from $L^\beta$ to $L^p$., and where $T_\lambda^{N-2k}$ acts on $L^\beta$. By Lemma \ref{action de T sur L^beta},

$$||T_\lambda^k||_{p,p}\leq C(1-\mu)^k,\,k\geq 2N,$$
which in turns implies that the series $\sum_{n\geq 0} T_\lambda^n$ converges in $\mathcal{L}(L^p,L^p)$ for all $1\leq p\leq\infty$, uniformly with respect to $\lambda\geq0$. This concludes the proof of Theorem \ref{resolvante du Laplacien de Hodge}.
\cqfd

\textit{Proof of Lemma \ref{action de T sur L^beta}:}\\\\ 
We write :

$$T_\lambda=H_\lambda^{-1/2}[H_\lambda^{-1/2}\mathcal{R}_-H_\lambda^{-1/2}]H_\lambda^{1/2},$$
and we define $A_\lambda:=H_\lambda^{-1/2}\mathcal{R}_-H_\lambda^{-1/2}$. Let us define the Hilbert space $H_{0,\lambda}^1$ to be the closure of $C_0^\infty(E)$ for the norm:

$$\omega\mapsto\left(\int_M|\nabla\omega|^2+\lambda|\omega|^2\right)^{1/2}=Q_\lambda(\omega)^{1/2},$$
where $Q_\lambda$ is the quadratic form associated to the self-adjoint operator $H_\lambda$. If $\lambda>0$, it is the space $H_0^1\cap L^2=\mathcal{D}om(H^{1/2})$, but with a different norm. The choice of the norm is made so that $H_\lambda^{1/2} : H_{0,\lambda}^1\rightarrow L^2$ is an isometry. Since $A_\lambda : L^2\rightarrow L^2$, and given that $T_\lambda=H_\lambda^{-1/2}A_\lambda H_\lambda^{1/2}$, we deduce that :

$$T_\lambda : H_{0,\lambda}^1\rightarrow H_{0,\lambda}^1\hbox{ with }||T_\lambda||_{H_{0,\lambda}^1,H_{0,\lambda}^1}=||A_\lambda||_{2,2}.$$
But by the equivalence $1)\Leftrightarrow 3)$ in the Definition \ref{sous-criticite}, the existence of $\mu\in (0,1)$ such that $||A_\lambda||_{2,2}\leq 1-\mu$ is equivalent to:

$$\langle \mathcal{R}_-\omega,\omega\rangle\leq (1-\mu)\langle (H_\lambda)\omega,\omega\rangle,\,\forall \omega\in C_0^\infty(\Gamma(E)).$$
Since $\langle (H+\lambda)\omega,\omega\rangle=\langle H\omega,\omega\rangle+\lambda||\omega||_2^2\geq\langle H\omega,\omega\rangle$, we obtain that the existence of some $\mu\in (0,1)$ such that for all $\lambda\geq0$, $||A_\lambda||_{2,2}\leq 1-\mu$ is equivalent to the strong positivity of $L$. Therefore $||T_\lambda||_{H_{0,\lambda}^1,H_{0,\lambda}^1}\leq (1-\mu)$. Moreover, by Theorem \ref{inegalites fonctionnelles},

$$H_\lambda^{-1/2} : L^\frac{2n}{n+2}\rightarrow L^2,$$
with norm bounded independantly of $\lambda\geq0$ (by Proposition \ref{independance of lambda}, and by Lemma \ref{action de T sur les L^p},

$$\mathcal{R}_- : L^\frac{2n}{n-2}\rightarrow L^\frac{2n}{n+2},$$
so that, using that $H_\lambda^{-1/2} : L^2\rightarrow H_{0,\lambda}^1$ is an isometry and that we can write $T_\lambda=H_\lambda^{-1/2}[H_\lambda^{-1/2}\mathcal{R}_-]$, we get that

$$T_\lambda : L^\beta\rightarrow H_{0,\lambda}^1,$$
is bounded with a bound of the norm independant of $\lambda\geq0$. Furthermore, $H_{0,\lambda}^1\hookrightarrow H_0^1$ is continuous of norm less than $1$, and the Sobolev inequality for $H$ says precisely that:

$$H_{0}^1\hookrightarrow L^\beta$$
continuously. Therefore, $H_{0,\lambda}^1\hookrightarrow L^\beta$ continuously, with a bound of the norm independant of $\lambda$. Then we write $T_\lambda^k= T_\lambda^{k-1}T_\lambda$, with

$$T_\lambda : L^\beta\rightarrow H_{0,\lambda}^1$$
and 

$$T_\lambda^{k-1} : H_{0,\lambda}^1\rightarrow H_{0,\lambda}^1\hookrightarrow L^\beta,$$
so that we get:

$$||T_\lambda^k||_{\beta,\beta}\leq C(1-\mu)^k.$$
\cqfd
As a byproduct of the proof (more precisely, of Proposition \ref{action de T^N sur les L^p} and Lemma \ref{action de T sur L^beta}), we get:

\begin{cor}\label{decomposition de la resolvante du Laplacien de Hodge}
 
$$(L+\lambda)^{-1}=(I-T_\lambda)^{-1}H_\lambda^{-1},$$
with $(I-T_\lambda)^{-1} : L^p(E)\rightarrow L^p(E)$ bounded with a bound of the norm independant of $\lambda$, for all $1\leq p \leq\infty$.

\end{cor}
We could hope to deduce from Theorem \ref{resolvante du Laplacien de Hodge} that $e^{-tL}$ is uniformly bounded on all the $L^p$ spaces, by an argument similar to the Hille-Yosida Theorem. In particular, the Hille-Yosida-Phillips Theorem tells us that the bound 

$$||(L+\lambda)^{-k}||\leq \frac{C}{\lambda^k},\,\forall k\in \mathbb{N},$$
with $C$ independant of $\lambda$ and $k$, is necessary and sufficient to obtain $e^{-tL}$ uniformly bounded. The issue here is that applying Theorem \ref{resolvante du Laplacien de Hodge} directly yields:

$$||(L+\lambda)^{-k}||\leq \frac{C^k}{\lambda^k},\,\forall k\in \mathbb{N},$$
i.e. the constant is not independant of $k$. In fact, applying the method of Theorem \ref{resolvante du Laplacien de Hodge} in a less naive way would in fact yield: 

$$||(L+\lambda)^{-k}||\leq \frac{Ck}{\lambda^k},\,\forall k\in \mathbb{N},$$ 
i.e. the growth of the constant is linear in $k$ and not exponential. Still this is not enough.\\
We will use an idea of Sikora to overcome this problem: it is shown in \cite{Sikora} that a Gaussian estimate for $e^{-tL}$ can be obtained using suitable on-diagonal estimates. Therefore, our goal will be to get these on-diagonal estimates for $e^{-tL}$, that is estimates for $||e^{-tL}||_{2,\infty}$, and in order to get these we can try, following Sikora, to get estimates on $||(L+\lambda)^{-k}||_{2,\infty}$. The point is that the bound needed on $||(L+\lambda)^{-k}||_{2,\infty}$ need not be independant of $k$, so Theorem \ref{resolvante du Laplacien de Hodge} should be enough to prove it! We follow this path in the next section.  

\begin{rem}
 
Of course, at the end, since we succeed in proving the Gaussian estimate for $e^{-tL}$, we obtain that $e^{-tL}$ is uniformly bounded on all the $L^p$ spaces.

\end{rem}

\subsection{On-diagonal upper bounds}

The next Proposition is a slight generalisation of Sikora's ideas:

\begin{prop}\label{estimee diagonale pour un operateur auto-adjoint}

Let $X$ be a measurable metric space. Let $L$ be a self-adjoint, positive unbounded operator on $L^2(X)$ , and let $1< p< \infty$. Assume that the semigroup $e^{-tL}$ is analytic bounded on $L^p(X)$ (it is necessarily the case if $p=2$). The following statements are equivalent:

\begin{enumerate}
 \item There exists a constant $C$ such that for all $t>0$,
 
 $$||e^{-tL}||_{p,\infty}\leq \frac{C}{t^{n/2p}}.$$
 
 \item For an (for all) $\alpha>n/2p$, there exists a constant $C(p,\alpha)$ such that 
 
 $$||(L+\lambda)^{-\alpha}||_{p,\infty}\leq C(p,\alpha)\lambda^{-\alpha+n/2p},\,\forall \lambda>0.$$ 
 
\end{enumerate}

\end{prop}
\pr First, notice that 

$$||(L+\lambda)^{-\alpha}||_{p,\infty}\leq C(p,\alpha)\lambda^{-\alpha+n/2p},\,\forall \lambda>0$$ 
can be rewritten as 

$$||(I+tL)^{-\alpha}||_{p,\infty}\leq C(p,\alpha)t^{-n/2p},\,\forall t>0.$$

\noindent $2)\Rightarrow 1)$: since $e^{-tL}$ is analytic bounded on $L^p$, by Proposition \ref{analyticity} there is a constant $C$ such that:

$$||(I+tL)^\alpha e^{-tL}||_{p,p}\leq C,\,\forall t>0.$$
We then write $e^{-tL}=(I+tL)^{-\alpha}\big((I+tL)^\alpha e^{-tL}\big)$ to obtain the result.

\vskip5mm

\noindent$1)\Rightarrow 2)$: we have 

$$(L+\lambda)^{-\alpha}=\frac{1}{\Gamma(\alpha+1)}\int_0^\infty e^{-\lambda t}e^{-tL}t^{\alpha-1}dt,$$
so that 

$$||(L+\lambda)^{-\alpha}||_{p,\infty}\leq\frac{1}{\Gamma(\alpha+1)}\int_0^\infty e^{-\lambda t}||e^{-tL}||_{p,\infty}t^{\alpha-1}dt.$$
Using the hypothesis, we obtain:

$$||(L+\lambda)^{-\alpha}||_{p,\infty}\leq\frac{1}{\Gamma(\alpha+1)}\int_0^\infty e^{-\lambda t}t^{\alpha-n/2p-1}dt=\frac{1}{\Gamma(\alpha+1)}\lambda^{-\alpha+n/2p}\int_0^\infty e^{-u}u^{\alpha-n/2p-1}du.$$
Since $\alpha-n/2p>0$, the integral $\int_0^\infty e^{-u}u^{\alpha-n/2p-1}du$ converges, hence the result.

\cqfd
We will use both sides of the equivalence. First, we apply this to $H$ (which, by Corollary \ref{e^-tH semi-groupe de contraction}, satisfies $||e^{-tH}||_{p,\infty}\leq\frac{C}{t^{n/2p}}$ and which is analytic bounded on $L^p$ for $1<p<\infty$ by Corollary \ref{e^-tH semi-groupe de contraction}), to get:

\begin{cor}\label{estimee diagonale de la resolvante de H}
 
For all $1\leq p\leq\infty$ and $\alpha>n/2p$, there exists a constant $C(p,\alpha)$ such that 
 
$$||H_\lambda^{-\alpha}||_{p,\infty}\leq C(p,\alpha)\lambda^{-\alpha+n/2p},\,\forall \lambda>0.$$ 

\end{cor}
We now use the other side of the equivalence in Proposition \ref{estimee diagonale pour un operateur auto-adjoint} (i.e. a bound on the resolvent implies a bound on the semigroup) to prove the following Theorem, which is our main result in this section:

\begin{thm}\label{estimee diagonale pour le Laplacien de Hodge}
 
Let $(M,g)$ be an complete Riemannian manifold which satisfies a Sobolev inequality of dimension $n$ \eqref{sobolev}, and assume that $\mathcal{R}_-$ is in $L^{\frac{n}{2}\pm\varepsilon}$ for some $\varepsilon>0$. We also assume that $L:=H-\mathcal{R}_-=\nabla^*\nabla+\mathcal{R}_+-\mathcal{R}_-$, acting on the sections of a fibre bundle $E\rightarrow M$, is strongly positive. Then we have the following on-diagonal estimate: there is a constant $C$ such that

$$||e^{-tL}||_{2,\infty}\leq \frac{C}{t^{n/4}},\,\forall t>0.$$

\end{thm}
\pr In this proof, we write $L^q$ for $L^q(E)$. By Proposition \ref{estimee diagonale pour un operateur auto-adjoint}, it is enough to prove the estimate:

\begin{equation}\label{estimee diagonale}
 ||(L+\lambda)^{-N}||_{2,\infty}\leq \frac{C_N}{\lambda^{N-n/4}},\,\forall\lambda>0,
\end{equation}

for an $N>n/4$. We use the fact that for all $1\leq p\leq\infty$, we have $(L+\lambda)^{-1}=(I-T_\lambda)^{-1}H_\lambda^{-1}$ on $L^p$, where $(I-T_\lambda)^{-1}$ is bounded on all the $L^p$ spaces, with a bound for the norm independant of $\lambda\geq0$ (c.f. Corollary \ref{decomposition de la resolvante du Laplacien de Hodge}). Let $k=\lfloor n/4\rfloor=\lfloor \frac{1}{2}/\frac{2}{n}\rfloor$. We will show the estimate (\ref{estimee diagonale}) for $N=k+1$.\\\\
\underline{First case:} $\frac{n}{4}\notin \mathbb{N}$\\
We want to show the estimate $||(L+\lambda)^{-k-1}||_{2,\infty}\leq \frac{C}{\lambda^{(k+1)-n/4}},\,\forall\lambda>0.$ Define $p>\frac{n}{2}$ by:

$$\frac{1}{p}=\frac{1}{2}-k\frac{2}{n}.$$
By the mapping property of Theorem \ref{inegalites fonctionnelles}, 

$$H_\lambda^{-1} : L^r\longrightarrow L^s,\,\frac{1}{s}=\frac{1}{r}-\frac{2}{n},\,\forall r<\frac{n}{2},$$
with a norm bounded independantly of $\lambda$. Using the fact that $(I-T_\lambda)^{-1}$ is bounded on all the $L^p$ spaces, with a bound for the norm independant of $\lambda\geq0$, we get that

$$(L+\lambda)^{-k} : L^2\longrightarrow L^p$$
is bounded uniformly in $\lambda\geq0$. Since $\frac{n}{2p}<1$, we have by Corollary \ref{estimee diagonale de la resolvante de H}:

$$H_\lambda^{-1} : L^p\longrightarrow L^\infty,$$
with 
$$||H_\lambda^{-1}||_{p,\infty}\leq C\lambda^{-1+\frac{n}{2p}},$$ 
so that:

$$||(L+\lambda)^{-k-1}||_{2,\infty}\leq C(k)\lambda^{-1+\frac{n}{2p}}=\frac{C(k)}{\lambda^{k+1-n/4}},$$
which is what we need.\\\\
\underline{Second case:} $\frac{n}{4}\in \mathbb{N}$ hence $k=\frac{n}{4}$. We write $H_\lambda^{-1}=H_\lambda^{-\alpha}H^{-1+\alpha}$, where $\alpha\in(0,1)$. Then by Proposition \ref{resolvante de H}, $||H_\lambda^{-1+\alpha}||_{2,2}\leq \frac{1}{\lambda^{1-\alpha}}$, and 

$$H_\lambda^{-\alpha} : L^2\longrightarrow L^q,\,\frac{1}{q}=\frac{1}{2}-\alpha\frac{2}{n}$$
is bounded with a norm bounded independantly of $\lambda>0$. This time, we define $p>\frac{n}{2}$ by:

$$\frac{1}{p}=\frac{1}{2}-(k-1+\alpha)\frac{2}{n}.$$
We get:

$$||(L+\lambda)^{-k}||_{2,p}\leq ||(L+\lambda)^{-(k-1)}||_{q,p}||(I-T_\lambda)^{-1}||_{q,q}||H_\lambda^{-\alpha}||_{2,q}||H_\lambda^{-(1-\alpha)}||_{2,2}\leq
\frac{C}{\lambda^{1-\alpha}},$$
Therefore, using that $||H_\lambda^{-1}||_{p,\infty}\leq C\lambda^{-1+\frac{n}{2p}}$ and $||(I-T_\lambda)^{-1}||_{\infty,\infty}\leq C$ independant of $\lambda$, we obtain:

$$||(L+\lambda)^{-k-1}||_{2,\infty}\leq \frac{C_k}{\lambda^{(1-n/2p)+(1-\alpha)}}.$$
But $\frac{n}{2p}=\frac{n}{4}-(k-1+\alpha)$, which yields what we want.
\cqfd

\subsection{Pointwise estimates of the heat kernel on 1-forms}
Let us recall the following definition:

\begin{defin}

Let $X$ be a metric measured space, $E$ a Riemannian vector bundle over $X$, and $L$ a non-negative self-adjoint operator on $L^2(E)$. We say that $L$ satisfies the \textbf{finite propagation speed property} if for every $t>0$, the support of the kernel of $\cos(t\sqrt{L})$ is included in the set $\{(x,y)\in X\times X : d(x,y)\leq t\}.$

\end{defin}
A consequence of Sikora's work (Theorem 4 in \cite{Sikora}) is:

\begin{thm}\label{resultat de Sikora}
 
Let $X$ be a metric measured space whose measure is doubling, and $E$ a Riemannian vector boundle over $X$. If the following on-diagonal estimate holds:

$$||e^{-tL}||_{2,\infty}\leq\frac{C}{t^{n/4}},\,\forall t>0,$$
with $L$ non-negative self-adjoint operator on $L^2(E)$, satisfying the finite speed propagation, then there is a Gaussian-type estimate for $e^{-tL}$: for every $\delta>0$, there is a constant $C$ such that

$$||K_{\exp(-tL)}(x,y)||\leq \frac{C}{t^{n/2}}\exp\left(-\frac{d^2(x,y)}{(4+\delta)t}\right),\,\forall x,y\in M,\,\forall t>0,$$
where $K_{\exp(-tL)}$ denotes the kernel of $e^{-tL}$.

\end{thm}
It is shown the following fact in the appendix of \cite{Ma-Marinescu}, p.388-389:

\begin{prop}

Let $M$ be a complete Riemannian manifold, $E$ a Riemannian vector bundle over $M$, and $L$ an operator of the type:

$$L:=\nabla^\star\nabla+\mathcal{R},$$
such that $L$ is non-negative self-adjoint on $L^2(E)$. Then $L$ satisfies the finite propagation speed property.

\end{prop}
Therefore, we have shown the following result, consequence of Theorems \ref{resultat de Sikora} and \ref{estimee diagonale pour le Laplacien de Hodge} :

\begin{thm}\label{estimee e^-tL}

Let $M$ be a complete Riemannian manifold satisfying the Sobolev inequality of dimension $n$ \eqref{sobolev},  $E$ a Riemannian vector bundle over $M$, and $L$ an operator of the type:

$$L:=\nabla^\star\nabla+\mathcal{R}_+-\mathcal{R}_-,$$
such that $L$ is self-adjoint on $L^2(E)$. We assume that $\mathcal{R}_-\in L^{\frac{n}{2}\pm\varepsilon}$, for some $\varepsilon>0$, and that $L$ is strongly positive. Then for every $\delta>0$, there is a constant $C$ such that

$$||K_{\exp(-tL)}(x,y)||\leq \frac{C}{t^{n/2}}\exp\left(-\frac{d^2(x,y)}{(4+\delta)t}\right),\,\forall x,y\in M,\,\forall t>0,$$

\end{thm}
The bound that we obtain is not exactly what is usually called a Gaussian estimate for $e^{-tL}$; indeed, a Gaussian estimate for $e^{-tL}$ is a bound of the following type:

$$||K_{\exp(-tL)}(x,y)||\leq \frac{C}{V(x,t^{1/2})}\exp\left(-\frac{d^2(x,y)}{(4+\delta)t}\right),\,\forall x,y\in M,\,\forall t>0.$$
The problem comes from the term $V(x,t^{1/2})$, which may not behave like $t^{-n/2}$. Indeed, when $M$ satisfies a Sobolev inequality of dimension $n$, we only have the lower bound (proved in \cite{Carron2} and \cite{Akutagawa}):

$$V(x,R)\geq CR^n,\,\forall R>0,\forall x\in M,$$
which implies by the way that $n\geq dim(M)$. For example the Heisenberg group $\mathbb{H}_1$ is a manifold of dimension $3$ which satisfies a Sobolev inequality of dimension 4, but whose volume of geodesic balls satisfies:

$$V(x,R)\approx R^3 \hbox{ if }R\leq1,$$
and 

$$V(x,R)\approx R^4 \hbox{ if }R\geq1.$$

\begin{defin}
 
We say that $M$ satisfies a \textbf{relative Faber-Krahn inequality} of exponent $n$ if there is a constant $C$ such that for every $x\in M$ and $R>0$, and every non-empty subset $\Omega\subset B(x,R)$,

$$\lambda_1(\Omega)\geq \frac{C}{R^2}\left(\frac{|B(x,R)|}{|\Omega|}\right)^{2/n},$$
where $\lambda_1(\Omega)$ is the first eigenvalue of $\Delta$ on $\Omega$ with Dirichlet boundary conditions. 

\end{defin}
It is proved in \cite{Grigor'yan2} that the relative Faber-Krahn inequality is equivalent to the volume doubling property \eqref{doubling} together with the Gaussian upper bound of the heat kernel \eqref{gaussian}:

$$p_t(x,y)\leq \frac{C}{V(x,\sqrt{t})}\exp\left(-\frac{d^2(x,y)}{(4+\delta)t}\right),\,\forall t>0,\forall x,y\in M.$$
We have the following property, which is not new but whose proof is given for the reader's convenience:

\begin{prop}\label{Faber-Krahn et Sobolev}

Let $M$ be a complete Riemannian manifold of dimension $m$, which satisfies a Sobolev inequality of dimension $n$ \eqref{sobolev}, and whose Ricci curvature is bounded from below. If the volume of big balls is Euclidean of dimension $n$:

$$V(x,R)\simeq C R^n,\,\forall R\geq1,$$
 then $M$ satisfies a relative Faber-Krahn inequality of exponent $n$.

\end{prop}
\pr Let us explain first why the relative Faber-Krahn inequality holds for balls of small radius. Saloff-Coste has shown in  \cite{Saloff-Coste2} the following Sobolev inequality: if the Ricci curvature of $M$ is bounded from below by$-K\leq0$, then for all ball $B$ of radius $R$,

\begin{equation}\label{sobolev petites boules}
||f||_{\frac{2n}{n-2}}^2\leq e^{C(1+\sqrt{K}R)}\frac{R^2}{V(R)^{2/n}}||df||_2+e^{C(1+\sqrt{K}R)}\frac{1}{V(R)^{2/n}}||f||_2,\,\forall f\in C_0^\infty(B)
\end{equation}
For balls of radius smaller than $1$, \eqref{sobolev petites boules} rewrites

$$||f||_{\frac{2n}{n-2}}^2\leq C\frac{R^2}{V(R)^{2/n}}||df||_2+C\frac{1}{V(R)^{2/n}}||f||_2,\,\forall f\in C_0^\infty(B).$$
Moreover, for balls of radius smaller than $1$, we have the following inequality for the first eigenvalue af the Laplacian with Dirichlet boundary conditions, consequence of Cheng's comparison Theorem:

$$\lambda_1(B)\leq C R^2,$$
therefore we obtain that for every ball of radius $R\leq 1$,

$$||f||_{\frac{2n}{n-2}}^2\leq C\frac{R^2}{V(R)^{2/n}}||df||_2,\,\forall f\in C_0^\infty(B).$$
From the work of Carron \cite{Carron2}, this is equivalent to the relative Faber-Krahn inequality for balls of radius smaller than $1$.\\\\
For balls of radius greater than $1$: again, according to \cite{Carron2}, since $M$ satisfies a Sobolev inequaity of dimension $n$ \eqref{sobolev}, $M$ satisfies a Faber-Krahn inequality of exponent $n$, that is for every open set $\Omega\subset M$,

$$\lambda_1(\Omega)\geq \frac{C}{|\Omega|^{\frac{2}{n}}}.$$
If $\Omega\subset B(x,R)$ with $R\geq1$, we have, using the hypothesis that the volume of balls of radius greater than $1$ is Euclidean of dimension $n$:

$$\lambda_1(\Omega)\geq \frac{C}{|\Omega|^{\frac{2}{n}}}\geq \frac{C}{R^2}\left(\frac{|B(x,R)|}{|\Omega|}\right)^{\frac{2}{n}}.$$
\cqfd

\begin{ex}
 
The Heisenberg group $\mathbb{H}_1$ satisfies a relative Faber-Krahn inequality of exponent $4$; in fact, it even satisfies the scaled Poincar\'{e} inequalities and the Doubling Property, which is equivalent (by the work of Grigor'yan \cite{Grigor'yan3} and Saloff-Coste \cite{Saloff-Coste}) to the conjunction of a Gaussian upper and lower bound for the heat kernel.\\\\
Every manifold with $Ric\geq0$ (or more generally, with $Ric\geq0$ outside a compact set, finite first Betti number, only one end, and satisfying a condition called (RCA), see \cite{Grigor'yan-Saloff-Coste}) satisfies the scaled Poincar\'{e} inequalities, and thus a relative Faber-Krahn inequality of exponent $dim(M)$.

\end{ex}
Taking into account what we have obtained in Theorem \ref{estimee diagonale pour le Laplacien de Hodge}, we get the result of Theorem \ref{gaussian estimate}, as anounced to in the introduction:

\begin{thm}\label{estimee gaussienne pour le Laplacien de Hodge}
 
Let $(M,g)$ be a complete Riemannian manifold which satisfies the Sobolev inequality of dimension $n$ \eqref{sobolev}, and $E$ a Riemannian vector bundle over $M$. Let $L$ a generalised Schr\"{o}dinger operator:

$$L:=\nabla^\star\nabla+\mathcal{R}_+-\mathcal{R}_-,$$
acting on the sections of $E$. We assume that $\mathcal{R}_-$ is in $L^{\frac{n}{2}-\varepsilon}\cap L^\infty$ for some $\varepsilon>0$, and that $L$ is strongly positive. We also assume that the volume of big balls is Euclidean of dimension $n$:

$$V(x,R)\simeq C R^n,\,\forall R\geq1.$$ 
Then the Gaussian estimate holds for $e^{-tL}$: for every $\delta>0$, there is a constant $C$ such that

$$||K_{\exp(-tL)}(x,y)||\leq \frac{C}{V(x,t^{1/2})}\exp\left(-\frac{d^2(x,y)}{(4+\delta)t}\right),\,\forall (x,y)\in M\times M,\,\forall t>0.$$

\end{thm}
\pr By Theorem \ref{estimee e^-tL},

$$||K_{\exp(-tL)}(x,y)||\leq \frac{C}{V(x,t^{1/2})}\exp\left(-\frac{d^2(x,y)}{(4+\delta)t}\right),\,\forall (x,y)\in M\times M,\,\forall t\geq1.$$
Since $M$ satisfies a relative Faber-Krahn inequality,

$${p_t}(x,y)\leq \frac{C}{V(x,t^{1/2})}\exp\left(-\frac{d^2(x,y)}{(4+\delta)t}\right),\,\forall (x,y)\in M\times M,\,\forall t>0.$$
But since $\mathcal{R}_-$ is bounded from below, this implies the Gaussian estimate for $e^{-tL}$ for small times:

$$||K_{\exp(-tL)}(x,y)||\leq \frac{C}{V(x,t^{1/2})}\exp\left(-\frac{d^2(x,y)}{(4+\delta)t}\right),\,\forall (x,y)\in M\times M,\,\forall t\leq1.$$
Indeed, this comes from the fact that we have the domination (proved in \cite{H-S-U1}) : 

$$||K_{\exp(-tL)}(x,y)||\leq e^{-t(\Delta-C)}(x,y)$$
if $\mathcal{R}_-\leq C$.
\cqfd
We will see in Proposition \ref{formes harmoniques Lp} that in fact, under the assumptions of Theorem \ref{estimee gaussienne pour le Laplacien de Hodge}, 

$$Ker_{L^2}(L)=Ker_{H_0^1}(L).$$
Using this and the definition of strong positivity, we get the result anounced in Theorem \ref{gaussian estimate}:

\begin{cor}\label{estimee gaussienne sans cohomologie}

Let $(M,g)$ be a complete Riemannian manifold which satisfies the Sobolev inequality of dimension $n$ \eqref{sobolev}, and $E$ a Riemannian vector bundle over $M$. Let $L$ be a generalised Schr\"{o}dinger operator:

$$L:=\nabla^\star\nabla+\mathcal{R}_+-\mathcal{R}_-,$$
acting on the sections of $E$, such that $L$ is non-negative on $L^2(E)$. We assume that $\mathcal{R}_-$ is in $L^{\frac{n}{2}-\varepsilon}\cap L^\infty$ for some $\varepsilon>0$, that the volume of big balls is Euclidean of dimension $n$:

$$V(x,R)\simeq C R^n,\,\forall R\geq1,$$ 
and that

$$Ker_{L^2}(L)=\{0\}.$$
Then the Gaussian estimate holds for $e^{-tL}$: for every $\delta>0$, there is a constant $C$ such that

$$||K_{\exp(-tL)}(x,y)||\leq \frac{C}{V(x,t^{1/2})}\exp\left(-\frac{d^2(x,y)}{(4+\delta)t}\right),\,\forall x,y\in M,\,\forall t>0.$$

\end{cor}

\section{Applications}

The Gaussian estimate on the heat kernel on 1-forms has a certain number of consequences, which we decribe now.

\subsection{Estimates on the gradient of the heat kernel on functions and scaled Poincar\'{e} inequalities}

We recall a classical definition:

\begin{defin}
We say that $M$ satisfies the \textbf{scaled Poincar\'{e} inequalities} if there exists a constant $C$ such that for every ball $B=B(x,r)$ and for every function $f$ with $f,\,\nabla f$ locally square integrable, 

$$\int_B|f-f_B|^2\leq C r^2\int_B|\nabla f|^2.$$

\end{defin}
Coulhon and Duong (p. 1728-1751 in \cite{Coulhon-Duong1}) have noticed that a Gaussian estimate on the heat kernel on 1-forms --in fact, a Gaussian estimate on the heat kernel on exact 1-forms is enough-- leads to the following estimate for the gradient of the heat kernel on functions:

$$|\nabla_x p_t(x,y)|\leq \frac{C}{\sqrt{t}V(x,\sqrt{t})}\exp\left(-\frac{d^2(x,y)}{(4+\delta)t}\right),\,\forall t>0,\,\forall x,y\in M,$$
which, when the on-diagonal Gaussian lower bound for the heat kernel on functions $p_t(x,x)\geq \frac{C}{V(x,\sqrt{t})}$ and the volume doubling property \eqref{doubling} hold, yields the Gaussian lower bound for the heat kernel on functions:

$$p_t(x,y)\geq \frac{C}{V(x,\sqrt{t})}\exp\left(-\frac{d^2(x,y)}{(4+\delta)t}\right),\,\forall t>0,\,\forall x,y\in M.$$
In addition, if $M$ satisfies a Sobolev inequality of dimension $n$ \eqref{sobolev} and if the volume growth of $M$ is compatible with the Sobolev dimension, we know from Proposition \ref{Faber-Krahn et Sobolev} that $M$ satisfies a relative Faber-Krahn inequality of exponent $n$, and this implies by the work of Grigor'yan \cite{Grigor'yan2} that we have the corresponding upper-bound for the heat kernel on functions:

$$p_t(x,y)\leq \frac{C}{V(x,\sqrt{t})}\exp\left(-\frac{d^2(x,y)}{(4+\delta)t}\right),\,\forall t>0,\,\forall x,y\in M$$
But we know from the work of Saloff-Coste and Grigor'yan \cite{Saloff-Coste}, \cite{Grigor'yan3} that the two-sided Gaussian estimates for the heat kernel on functions are equivalent to the conjonction of the scaled Poincar\'{e} inequalities and the volume doubling property \eqref{doubling}.\\
Thus we have proved the following theorem, which extends similar results for manifolds with non-negative Ricci curvature:

\begin{thm}\label{estimee du gradient du noyau de la chaleur}
 
Let $(M,g)$ be a complete Riemannian manifold which satisfies a Sobolev inequality of dimension $n$ \eqref{sobolev}, and whose negative part of the Ricci tensor $Ric_-$ is in $L^{\frac{n}{2}\pm\varepsilon}$ for some $\varepsilon>0$. We assume that there is no non-zero $L^2$ harmonic 1-form on $M$, that the volume of big balls is Euclidean of dimension $n$:

$$V(x,R)\simeq C R^n,\,\forall R\geq1,$$
and that the Ricci curvature is bounded from below. Then we have the following estimates on the heat kernel on functions:

$$|\nabla_x p_t(x,y)|\leq \frac{C}{\sqrt{t}V(x,\sqrt{t})}\exp\left(-\frac{d^2(x,y)}{(4+\delta)t}\right),\,\forall t>0,\,\forall x,y\in M,$$

$$\frac{c}{V(x,\sqrt{t})}\exp\left(-\frac{d^2(x,y)}{(4+\delta)t}\right)\leq p_t(x,y)\leq \frac{C}{V(x,\sqrt{t})}\exp\left(-\frac{d^2(x,y)}{(4+\delta)t}\right),\,\forall t>0,\,\forall x,y\in M,$$
and on $M$ the scaled Poincar\'{e} inequalities hold.

\end{thm}

\subsection{Boundedness of the Riesz transform}

In \cite{Sikora}, as an extension of classical Calder\'{o}n-Zygmund theory, A. Sikora shows that when a Gaussian estimate holds for a semigroup $e^{-tH}$, then for every local operator $A$ such that $AL^{-\alpha}$ is bounded on $L^2$ (with $\alpha>0$), the operator $AL^{-\alpha}$ is bounded on $L^p$ for all $1<p\leq 2$. Given this, we obtain the following corollaries, which are consequences of Theorem 10 in \cite{Sikora} (or Theorem 5.5 in Coulhon-Duong \cite{Coulhon-Duong1}, or the main result of \cite{Auscher-Coulhon-Duong-Hofmann}), and of Theorem \ref{estimee gaussienne pour le Laplacien de Hodge}: first, we have the result announced in Theorem \ref{Riesz 1<p<infty}

\begin{cor}\label{transformee de Riesz bornee dans le cas sous-critique}
 
Let $(M,g)$ be a complete Riemannian manifold which satisfies the Sobolev inequality of dimension $n$ \eqref{sobolev}, and whose negative part of the Ricci tensor $Ric_-$ is in $L^{\frac{n}{2}\pm\varepsilon}$ for some $\varepsilon>0$. We also assume that there is no non-zero $L^2$ harmonic 1-form on $M$, that the volume of big balls is Euclidean of dimension $n$

$$V(x,R)\simeq CR^n,\,\forall R\geq1,$$
and that the Ricci curvature is bounded from below. Then the Riesz transform $d\Delta^{-1/2}$ on $M$ is bounded on $L^p$, for all $1<p<\infty$.

\end{cor}
Then, we also have the following result which will be used in the section about the boundedness of the Riesz transform on $L^p$ for $1<p<n$:

\begin{cor}\label{Riesz potential}

Let $(M,g)$ be a complete Riemannian manifold which satisfies the Sobolev inequality of dimension $n$ \eqref{sobolev}, and whose negative part of the Ricci tensor $Ric_-$ is in $L^{\frac{n}{2}\pm\varepsilon}$ for an $\varepsilon>0$. We assume that the volume of big balls is Euclidean of dimension $n$

$$V(x,R)\simeq CR^n,\,\forall R\geq1,$$
and that the Ricci curvature is bounded from below. If $V$ is a non-negative potential such that $L:=\vec{\Delta}+V$ is strongly positive, then $d^\star\left(\vec{\Delta}+V\right)^{-1/2}$ is bounded on $L^p$ for every $1<p\leq 2$.

\end{cor}

\subsection{$L^p$ reduced cohomology}

\begin{defin}

The first space of $L^p$ \textbf{reduced cohomology}, denoted $H_p^1(M)$, is the quotient of $\{\alpha\in L^p(\Lambda^1T^\star M) : d\alpha=0\}$ by the closure in $L^p$ of the space of exact forms $dC_0^\infty(M)$. 

\end{defin}
Let us recall the following result, from \cite{Carron-Coulhon-Hassell}:

\begin{prop}\label{hodge cohomologie}

Let $p\geq 2$. Let $M$ be a complete non-compact Riemannian manifold, satisfying a Sobolev inequality of dimension $n$ \eqref{sobolev}, and whose volume of big balls is Euclidean of dimension $n$

$$V(x,R)\simeq CR^n,\,\forall R\geq1.$$ 
Assume that on $M$ the Riesz transform is bounded on $L^p$. Then $H_p^1(M)$ has the following interpretation:

\begin{equation}
H^1_p(M)\simeq \{\omega\in L^p(\Lambda^1 T^\star M) : d\omega=d^\star\omega=0\}
\end{equation}
Also, $\mathcal{H}^1(M)$, the space of $L^2$ harmonic forms, injects into $H^1_p(M)$.

\end{prop}
In particular, this implies that $H_p^1(M)$ is a space of harmonic forms:

$$H^1_p(M)\subset \{\omega\in L^p(\Lambda^1 T^\star M) : \vec{\Delta}\omega=0\}.$$
And furthermore,

\begin{prop}\label{formes harmoniques Lp}

Let $M$ be a complete, non-compact manifold, satisfying a Sobolev inequality of dimension $n$ \eqref{sobolev}, and whose volume of big balls is Euclidean of dimension $n$

$$V(x,R)\simeq CR^n,\,\forall R\geq1.$$ 
Let $E$ be a Riemannian vector bundle over $M$, endowed with a compatible connection $\nabla$, and $L$ a generalised Schr\"{o}dinger operator:

$$L=\nabla^\star\nabla+\mathcal{R},$$
acting on the sections of $E$. We assume that $\mathcal{R}_-$ is in $L^{\frac{n}{2}\pm\varepsilon}$ for some $\varepsilon>0$. Let $p\geq2$, then every section $\omega$ of $E$, which lies in $L^p$ and satisfies $L\omega=0$, is in $L^1\cap L^\infty$ (so in particular is in $L^2$). 

\end{prop}
In particular, for $L$ the Hodge Laplacian on 1-forms:

\begin{cor}

Let $p\geq 2$. Let $M$ be a complete non-compact Riemannian manifold, satisfying a Sobolev inequality of dimension $n$ \eqref{sobolev}, and whose volume of big balls is Euclidean of dimension $n$

$$V(x,R)\simeq CR^n,\,\forall R\geq1.$$
Assume that the negative part of the Ricci curvature is in $L^{\frac{n}{2}\pm\varepsilon}$ for some $\varepsilon>0$, and that the Riesz transform on $M$ is bounded on $L^p$. Then

$$H_1^p(M)=\mathcal{H}^1(M).$$

\end{cor}

\begin{rem}

This improves a result of Carron \cite{Carron3}, according to which if $M$ satisfies a Sobolev inequality of dimension $n$, and if the negative part of the Ricci curvature is in $L^{n/2}$ for $n>4$, then every $L^p$ harmonic form, for $p=\frac{2n}{n-2}$, is in $L^2$.

\end{rem}
\pr Let $\omega$ a section of $E$ in $L^p$, such that $L\omega=0$: that is

$$(\bar{\Delta}+\mathcal{R}_+)\omega-\mathcal{R}_-\omega=0.$$
Let $H:=\bar{\Delta}+\mathcal{R}_+$.

\begin{lem}

The following formula holds in $L^p$:

\begin{equation}\label{hodge Lp}
\omega=-H^{-1}\mathcal{R}_-\omega
\end{equation}

\end{lem}
\pr Let $\eta:=-H^{-1}\mathcal{R}_-\omega$, then $\eta\in L^p$ since according to Lemma \ref{action de T sur les L^p}, we have $H^{-1}\mathcal{R}_-\in \mathcal{L}(L^p,L^p)$, and furthermore

$$H(\omega-\eta)=0.$$
By Kato's inequality,

$$\Delta|\omega-\eta|\leq0,$$
i.e. $|\omega-\eta|$ is sub-harmonic. But according to Yau \cite{Yau}, there is no non-constant $L^p$ non-negative sub-harmonic functions on a complete manifold. $M$ being of infinite volume by the volume growth assumption, the only constant function in $L^p$ is the zero function. So we deduce that $\omega=\eta$.
\cqfd
\textit{End of the proof of Proposition \ref{formes harmoniques Lp} :} If we let $T:=H^{-1}\mathcal{R}_-$, then by Proposition \ref{action de T^N sur les L^p} there is an $N\in \mathbb{N}$ such that $T^N\omega\in L^1\cap L^\infty$. But by Lemma \ref{hodge Lp}, 

$$\omega=T^N\omega,$$
which shows that $\omega\in L^1\cap L^\infty$.
\cqfd

\section{Boundedness of the Riesz transform in the range $1<p<n$}

As announced in the introduction, we now remove the hypothesis that the space of $L^2$-harmonic 1-forms is reduced to $\{0\}$. We are mainly inspired by the perturbation technique developped by Carron in \cite{Carron4}. This section is devoted to the proof of Theorem \ref{Riesz 1,n}, which we recall here below:

\begin{thm}\label{transformee de Riesz pour 1<p<n}

Let $(M,g)$ be a complete Riemannian manifold which, for some $n>3$, satisfies a Sobolev inequality \eqref{sobolev}, and whose negative part of the Ricci tensor $Ric_-$ is in $L^{\frac{n}{2}\pm\epsilon}$ for an $\epsilon>0$. We also assume that the Ricci curvature is bounded from below, and that the volume of big balls is Euclidean of dimension $n$:

$$V(x,R)\simeq CR^n,\,\forall R\geq 1.$$
Then for every $1<p<n$, the Riesz transform is bounded on $L^p$ on $M$.

\end{thm}
We now give the proof of Theorem \ref{transformee de Riesz pour 1<p<n}, assuming two results (Theorem \ref{transformee de Riesz avec potentiel} and Theorem \ref{resultat de perturbation}) whose proof will be given in the next two subsections.\\\\
\textit{Proof of Theorem \ref{transformee de Riesz pour 1<p<n}:}\\\\

The hypotheses made imply by Proposition \ref{Faber-Krahn et Sobolev} that $M$ satisfies the relative Faber-Krahn inequality of exponent $n$, which is equivalent to the conjunction of the volume doubling property \eqref{doubling} and of the Gaussian upper-estimate \eqref{gaussian} on $p_t$, and we know by \cite{Coulhon-Duong2} that all this imply that the Riesz transform on $M$ is bounded on $L^p$ for all $1<p\leq2$. What we prove below is that the Riesz transform is bounded on $L^p$ for every $\frac{n}{n-1}<p<n$, which is thus enough to get the result.
As explained in the introduction, the proof is by a perturbation argument: using ideas of \cite{Carron4}, we will show in Theorem \ref{resultat de perturbation} that if $V\in C_0^\infty$ is non-negative, then $d(\Delta+V)^{-1/2}-d\Delta^{-1/2}$ is bounded on $L^p$ for $\frac{n}{n-1}<p<n$. Then we will prove in Theorem \ref{transformee de Riesz avec potentiel} that if $V$ is chosen such that $\vec{\Delta}+V$ is strongly positive, then $d(\Delta+V)^{-1/2}$ is bounded on $L^p$ for $\frac{n}{n-1}<p<n$. Finally, the following Lemma will conclude the proof of Theorem \ref{transformee de Riesz pour 1<p<n}:

\begin{lem}\label{ajout d'un potentiel positif}
 
Let $(M,g)$ be a complete Riemannian manifold which satisfies the Sobolev inequality of dimension $n$, and whose negative part of the Ricci tensor is in $L^{n/2}$. Then we can find a non-negative potential $V\in C_0^\infty$ such that $\vec{\Delta}+V$ is strongly positive.

\end{lem}
\textit{Proof of Lemma \ref{ajout d'un potentiel positif}:}\\\\

If we write $\vec{\Delta}+V=(\bar{\Delta}+W_+)-W_-=H-W_-$, and $A:=H^{-1/2}W_-H^{-1/2}$, then by the definition of strong positivity, $\vec{\Delta}+V$ is strongly positive if and only if $||A||_{2,2}<1$. Moreover, by Lemma (\ref{cas Ric dans L^n/2}), we have $||A||_{2,2}\leq C ||W_-||_{n/2}$, where $C$ is independant of the chosen potential $V\geq 0$. Therefore, it is enough to take $V$ such that $||(V-Ric_-)_-||_{\frac{n}{2}}<\frac{1}{C}$, which is possible since $Ric_-\in L^{n/2}$.
\cqfd
In the next two subsections, we prove the two results alluded to above: in the first one, we study the boundedness of $d(\Delta+V)^{-1/2}-d\Delta^{-1/2}$, and in the second one, the boundedness of $d(\Delta+V)^{-1/2}$.

\subsection{Boundedness of $d(\Delta+V)^{-1/2}-d\Delta^{-1/2}$}

Our aim here is to prove:

\begin{thm}\label{resultat de perturbation}

Let $(M,g)$ be a complete Riemannian manifold which satisfies a Sobolev inequality of dimension $n$ \eqref{sobolev}, and whose Ricci curvature is bounded from below. Let $V\in C_0^\infty$ be a non-negative potential. Then for every $\frac{n}{n-1}<p<n$, $d(\Delta+V)^{-1/2}-d\Delta^{-1/2}$ is bounded on $L^p$ on $M$.

\end{thm}
The proof is an adaptation of the proof in \cite{Carron4}. In order to adapt these ideas to the case of a Schr\"{o}dinger operator with non-negative potential, we will need some preliminary results. First, we recall an elliptic regularity result:

\begin{prop}\label{transformees de Riesz locales}
 
Let $V\in C_0^\infty$ be non-negative, and let $\Omega$ be a smooth, open, relatively compact subset. Let $\Delta_D$ be the Laplacian with Dirichlet conditions on $\Omega$.
Then the Riesz transforms $d(\Delta_D+V)^{-1/2}$ and $d\Delta_D^{-1/2}$ are bounded on $L^p$ for $1<p<\infty$. 

\end{prop}
We also recall the next Lemma and its proof from \cite{Carron4} :

\begin{lem}\label{lemme de Bakry}
 
Let $(M,g)$ be a complete Riemannian manifold with Ricci curvature bounded from below, and $V\in C_0^\infty$ be a non-negative potential. Then for all $1<p<\infty$, there is a constant $C$ such that  

$$||df||_p\leq C(||\Delta f||_p+||f||_p),\,\forall f\in C_0^\infty(M),$$
and 

$$||df||_p\leq C(||(\Delta+V) f||_p+||f||_p),\,\forall f\in C_0^\infty(M).$$

\end{lem}
\pr By Theorem 4.1 in \cite{Bakry}, the local Riesz transform is bounded on the $L^p$ for $1<p<\infty$, i.e. we have the following inequality for $a\geq0$ sufficently large:

$$||df||_p\leq C(||\Delta^{1/2} f||_p+a||f||_p),\,\forall f\in C_0^\infty(M).$$
We then use the fact that for all $1<p<\infty$, there exists a constant $C$ such that:

$$||\Delta^{1/2}f||_p\leq C\sqrt{||\Delta f||_p||f||_p}\leq \frac{C}{2} (||\Delta f||_p+||f||_p).$$
A proof of this inequality can be found in \cite{Coulhon-Russ-Tardivel}.\\
For the case with a potential, we have $||(\Delta+V)f||_p+a||f||_p\geq ||\Delta f||_p-||V||_\infty||f||_p+a||f||_p$. Taking $a>||V||_\infty$, we get the result.
\cqfd
\textit{Proof of Theorem \ref{resultat de perturbation}:}\\\\ 

Let $p\in (\frac{n}{n-1},n)$. We follow the proof of Carron in \cite{Carron4}. We define $L_0:=\Delta+V$, $L_1:=\Delta$; we take $K_1$ smooth, compact containing the support of $V$, and $K_2$, $K_3$ smooth, compact such that $K_1\subset\subset K_2\subset\subset K_3$. We also denote $\Omega:=M\setminus K_1$. Let $(\rho_0,\rho_1)$ a partition of unity such that supp $\rho_0\subset\Omega$ and supp $\rho_1\subset K_2$ (supp being for the support). We also take $\phi_0$ and $\phi_1$ to be $C^\infty$ non-negative functions such that supp $\phi_0\subset\Omega$, supp $\phi_1\subset K_3$ and such that $\phi_i\rho_i=\rho_i$. Moreover, we assume that $\phi_1 |_{K_2}=1$. \\
We define $H_0:=\Delta+V$ with Dirichlet boundary conditions on $K_3$, and $H_1:=\Delta$ with Dirichlet boundary conditions on $K_3$. Then, following Carron \cite{Carron4}, we construct parametrices for $e^{-t\sqrt{L_1}}$ and $e^{-t\sqrt{L_0}}$: the one for $e^{-t\sqrt{L_1}}$ is defined by

$$E_t^1(u):=\phi_1e^{-t\sqrt{H_1}}(\rho_1 u)+\phi_0e^{-t\sqrt{L_1}}(\rho_0 u),$$
and the one for $e^{-t\sqrt{L_0}}$ is defined by

$$E_t^0(u):=\phi_1e^{-t\sqrt{H_0}}(\rho_1 u)+\phi_0e^{-t\sqrt{L_1}}(\rho_0 u).$$
Let us note that for $e^{-t\sqrt{L_0}}$, we approximate by $e^{-t\sqrt{L_1}}$  outside the compact $K_3$, and not by $e^{-t\sqrt{L_0}}$. Let us also remark that $E_0^1(u)=E_0^0(u)=u$, as it should. We then have:

$$e^{-t\sqrt{L_i}}(u)=E_t^i(u)-G_i\big[(-\frac{\partial^2}{\partial t^2}+L_i)E_t^i(u)\big],$$
where $G_i$ is the Green operator on $\R_+\times M$ with Dirichlet boundary condition, associated to $-\frac{\partial^2}{\partial t^2}+L_i$.\\
Define 

$$(g_i(u))(x):=\int_0^\infty (G_i(S_t^i(u)))(t,x)dt.$$
Then, we have the following estimate whose proof is postponed to the end of this paragraph:

\begin{prop}\label{estimee derivee}
For any $\frac{n}{n-1}<p<n$, there is a constant $C$ (depending on $p$) such that for all $u\in L^p$,
$$||d(g_i(u))||_p\leq C||u||_p.$$
\end{prop}
Given this result, we can now finish the proof of Theorem \ref{resultat de perturbation}. We use the formula

$$L_i^{-1/2}=c \int_0^\infty e^{-t\sqrt{L_i}}dt,$$
to get:

$$L_i^{-1/2}u=\phi_1H_i^{-1/2}\rho_1 u+\phi_0 L_1^{-1/2}\rho_0 u-cg_i(u).$$
Therefore:

$$dL_1^{-1/2}u-dL_0^{-1/2}u=\big(d(\phi_1H_1^{-1/2}\rho_1 u)-d(\phi_1H_0^{-1/2}\rho_1 u)\big)+c\big(dg_0(u)-dg_1(u)\big).$$
(here is where we use the fact that we have taken for parametrices $e^{-t\sqrt{L_1}}$ for both operators outside a compact set). Write $d(\phi_1H_i^{-1/2}\rho_1 u)=(d\phi_1)H_i^{-1/2}\rho_1 u+\phi_1 dH_i^{-1/2}\rho_1 u$. $(d\phi_1)H_i^{-1/2}\rho_1$ has a smooth kernel with compact support, therefore is bounded on $L^p$. Applying Proposition \ref{transformees de Riesz locales}, we get that $\phi_1 dH_i^{-1/2}\rho_1$ is bounded on $L^p$. Also, Proposition \ref{estimee derivee} gives that $dg_0(u)$ and $dg_1(u)$ are bounded on $L^p$, therefore we get the result of Theorem \ref{resultat de perturbation}.
\cqfd
\textit{Proof of Proposition \ref{estimee derivee}:}\\\\
We compute:

$$\left(-\frac{\partial^2}{\partial t^2}+L_1\right)E_t^1(u)=[L_1,\phi_0]e^{-t\sqrt{L_1}}(\rho_0 u)+[L_1,\phi_1]e^{-t\sqrt{H_1}}(\rho_1 u),$$
and 

$$\left(-\frac{\partial^2}{\partial t^2}+L_0\right)E_t^0(u)=[L_0,\phi_1]e^{-t\sqrt{H_0}}(\rho_1 u)+[L_1,\phi_0]e^{-t\sqrt{L_1}}(\rho_0 u)+(L_0-L_1)\phi_0e^{-t\sqrt{L_1}}(\rho_0 u).$$
But $L_0-L_1=V$ is supported in $K_1$, therefore $(L_0-L_1)\phi_0e^{-t\sqrt{L_1}}(\rho_0 u)=0$. Moreover, we have $[\Delta+V,\phi_i]=[\Delta,\phi_i]$, therefore $[L_0,\phi_i]e^{-t\sqrt{H_0}}(\rho_i u)=(\Delta\phi_i)(e^{-t\sqrt{H_0}}(\rho_i u))-2\langle d\phi_i,\nabla e^{-tH_0}(\rho_i u)\rangle$. Define $S_t^i(u):=(-\frac{\partial^2}{\partial t^2}+L_i)E_t^i(u)$. We get:

$$S_t^1(u)=[\Delta,\phi_0]e^{-t\sqrt{L_1}}(\rho_0 u)+[\Delta,\phi_1]e^{-t\sqrt{H_1}}(\rho_1 u),$$
and

$$S_t^0(u)=[\Delta,\phi_0]e^{-t\sqrt{L_1}}(\rho_0 u)+[\Delta,\phi_1]e^{-t\sqrt{H_0}}(\rho_1 u).$$
Next,

\begin{lem}\label{terme d'erreur}
We have the following estimate for the error terms:

$$||S_t^i(u)||_1+||S_t^i(u)||_p\leq \frac{C}{(1+t)^{n/p}}||u||_p,\,\forall t>0.$$

\end{lem}
\pr Lemma 2.4 in \cite{Carron4} implies:

$$||[\Delta,\phi_0]e^{-t\sqrt{\Delta}}(\rho_0 u)||_1+||[\Delta,\phi_0]e^{-t\sqrt{\Delta}}(\rho_0 u)||_p\leq \frac{C}{(1+t)^{n/p}}.$$
Furthermore, if $f_1(u):=[\Delta,\phi_1]e^{-t\sqrt{H_1}}(\rho_1 u)=(\Delta\phi_1)e^{-t\sqrt{H_1}}(\rho_1 u)-2\langle d\phi_1,\nabla e^{-tH_1}(\rho_1 u)\rangle $,\\

and $f_0(u):=[\Delta,\phi_1]e^{-t\sqrt{H_0}}(\rho_1 u)=(\Delta\phi_1)e^{-t\sqrt{H_0}}(\rho_1 u)-2\langle d\phi_1,\nabla e^{-tH_0}(\rho_1 u)\rangle$, we have as in \cite{Carron4}:

$$||f_i(u)||_1+||f_i(u)||_p\leq \frac{C}{(1+t)^{n/p}}||u||_p,\,\forall t>0.$$
Indeed, if we denote $p_i^D(t,x,y)$ the heat kernel of $H_i$, then for $F_1$, $F_2$ disjoint compact subsets,

$$\lim_{t\rightarrow0} p_i^D(t,.,.)|_{F_1\times F_2} =0\hbox{ in }C^1$$
(cf \cite{Dodziuk} Lemma 3.2 and \cite{Ray-Singer}, Proposition 5.3). But by our hypotheses, the supports of $\rho_1$ and of $\Delta\phi_1$ are compact and disjoints, as are the ones of $\rho_1$ and $d\phi_1$. Therefore the kernels of the operators $S^i(t):=[\Delta,\phi_1]e^{-t\sqrt{H_i}}\rho_1$ are uniformly bounded as $t\rightarrow 0$. So we get:

$$||S^i(t)||_{p,\infty}\leq C,\,\forall t\in [0,1].$$
Now, the operators $H_i$ have a spectral gap, so $||e^{-t\sqrt{H_i}}||_{2,2}\leq e^{-ct}$, where $c>0$. If $v\in W^{1,2}(K_3)$ is a non-negative solution of $\frac{\partial v}{\partial t}+(\Delta+V)v=0$, then $\frac{\partial v}{\partial t}+\Delta v\leq 0$, and therefore by the parabolic maximum principle, $v$ attains its maximum on $\{t=0\}\cup \partial K_3$. If we take $v:=e^{-t(\Delta_D+V)}1$, which is zero on $\partial K_3$ for $t>0$, we get: 

$$\int_{K_3}p_i(t,x,y)dy\leq 1,\,\forall x\in K_3,$$
and therefore $||e^{-tH_i}||_{\infty,\infty}\leq 1$. By duality, it is true also on $L^1$, and by the subordination identity we have:
 
$$||e^{-t\sqrt{H_i}}||_{1,1}+||e^{-t\sqrt{H_i}}||_{\infty,\infty}\leq C.$$
Interpolating this with the $L^2$ bound, we get that

$$||e^{-t\sqrt{H_i}}||_{p,p}\leq Ce^{-ct},$$
for $1<p<\infty$, where the constants $C$ and $c$ depend on $p$. Then we write for $t\geq 1$:

$$||S^i(t)u||_\infty\leq ||[\Delta,\phi_1]e^{-\frac{1}{2}\sqrt{H_i}}||_{L^p\rightarrow L^\infty}||e^{-(t-1/2)\sqrt{H_i}}\rho_1 u||_{L^p}\leq Ce^{-ct}||u||_p.$$
Here we have used that the heat kernels $p_i^D(\frac{1}{2},.,.)$ are $C^\infty$.
\cqfd
\textit{End of the proof of Proposition \ref{estimee derivee}:}\\\\
To conclude the proof, we use an idea of G. Carron: the main argument of \cite{Carron4} together with the estimate of Lemma \ref{terme d'erreur} shows that 

$$||L_i(g_i(u))||_p+||g_i(u)||_p\leq C||u||_p.$$  
Applying Lemma \ref{lemme de Bakry}, we deduce that 

$$||d(g_i(u))||_p\leq C||u||_p.$$
\cqfd

\subsection{Boundedness of $d(\Delta+V)^{-1/2}$}

We now show:

\begin{thm}\label{transformee de Riesz avec potentiel}

Let $(M,g)$ be a complete Riemannian manifold which, for some $n>3$, satisfies the Sobolev inequality of dimension $n$ \eqref{sobolev}, and whose negative part of the Ricci tensor is in $L^{\frac{n}{2}-\varepsilon}\cap L^\infty$ for some $\varepsilon>0$. We also assume that the volume of big balls of $M$ is Euclidean of dimension $n$:

$$V(x,R)\simeq C R^n,\,\forall R\geq1.$$
Let $V\in C_0^\infty$ be non-negative, such that $\vec{\Delta}+V$ is strongly positive. Then the Riesz transform $d(\Delta+V)^{-1/2}$ is bounded on $L^p$ for every $1<p<n$.

\end{thm}
We first mention a preliminary result:

\begin{lem}\label{(Delta+V)^-1/2 d borne}
Under the assumptions of Theorem \ref{transformee de Riesz avec potentiel}, $(\vec{\Delta}+V)^{-1/2}d$ is bounded on $L^p$ for every $2\leq p<\infty$.

\end{lem}
\pr It is a direct consequence of Corollary \ref{Riesz potential}, by taking duals.
\cqfd
Let us now describe the idea of the proof of Theorem \ref{transformee de Riesz avec potentiel}: if the commutation relation $(\vec{\Delta}+V)^{-1/2}d=d(\Delta+V)^{-1/2}$ were true, we could conclude at once that  the Riesz transform with potential $d(\Delta+V)^{-1/2}$ is bounded on $L^p$ for all $2\leq p<\infty$. However, this is not true as soon as $V$ is not identically zero. Yet, we will prove Theorem \ref{transformee de Riesz avec potentiel} by proving that the error term $(\vec{\Delta}+V)^{-1/2}d-d(\Delta+V)^{-1/2}$ remains bounded on $L^p$ when $p<n$, and for this purpose we use once again the method of \cite{Carron4}.\\\\
\textit{Proof of Theorem \ref{transformee de Riesz avec potentiel}:}\\\\ First, let us note that we can restrict ourselves to the case $\frac{n}{n-1}<p<n$. Indeed, for $1<p<2$, since the hypotheses that we have made imply the Faber-Krahn inequality, and given the domination $e^{-t(\Delta+V)}\leq e^{-t\Delta}$, we have a Gaussian upper bound for $e^{-t(\Delta+V)}$. Thus the result of \cite{Coulhon-Duong2} shows that $d(\Delta+V)^{-1/2}$ is bounded on $L^p$ for every $1<p\leq 2$. So let $p\in(\frac{n}{n-1},n)$. We will use the following:

\begin{lem}\label{ultracontractivite pour le Laplacien de Hodge}
 
For $1\leq r\leq s\leq\infty$, we have the existence of a constant $C$ such that:

$$||e^{-t(\vec{\Delta}+V)}||_{L^r\rightarrow L^s}\leq\frac{C}{t^{\frac{n}{2}(\frac{1}{r}-\frac{1}{s})}}$$
and

$$||e^{-t\sqrt{\vec{\Delta}+V}}||_{L^r\rightarrow L^s}\leq\frac{C}{t^{n(\frac{1}{r}-\frac{1}{s})}}.$$

\end{lem}
We postpone the proof of Lemma \ref{ultracontractivite pour le Laplacien de Hodge} until the end of this section. Let $E$ be the vector bundle of basis $M\times \mathbb{R}_+$, whose fiber in $(t,p)$ is $\Lambda^1T_p^*M$. Let $G$ be the operator (the "Green operator'') acting on sections of $E$, whose kernel is given by 

$$G(\sigma,s,x,y)=\int_0^\infty \left[\frac{e^{-\frac{(\sigma-s)^2}{4t}}-e^{-\frac{(\sigma+s)^2}{4t}}}{\sqrt{4\pi t}}\right] \vec{p_t}^V(x,y)dt,$$
where $\vec{p_t}^V$ is the kernel of $e^{-t(\vec{\Delta}+V)}$. We can see that $G$ satisfies:

$$\left(-\frac{\partial^2}{\partial \sigma^2}+(\vec{\Delta}_x+V)\right)G=I,$$
and that $G(\sigma,s,x,y)$ is finite if $x\neq y$ and $\sigma\neq s$ (here we use the estimate $||\vec{p}_t^V(x,y)||\leq \frac{C}{t^{n/2}}$, given by Theorem \ref{estimee gaussienne pour le Laplacien de Hodge}). We want to prove, as in the proof of Theorem \ref{resultat de perturbation}, the following

\begin{prop}\label{parametrix}
The following formula holds for $u\in C^\infty_0$:

$$e^{-t\sqrt{\vec{\Delta}+V}}du=de^{-t\sqrt{\Delta+V}}u-G\left(\left(-\frac{\partial^2}{\partial t^2}+(\vec{\Delta}+V)\right)de^{-t\sqrt{\Delta+V}}u\right).$$
\end{prop}
The proof of Proposition \ref{parametrix} uses some estimates on the error term $G\left(\left(-\frac{\partial^2}{\partial t^2}+(\vec{\Delta}+V)\right)de^{-t\sqrt{\Delta+V}}u\right)$ to be proven below (Lemma \ref{majoration du terme d'erreur} and Lemma \ref{Green estimates}), and we postpone its proof for now. In order to estimate the error term, we begin with

\begin{lem}\label{majoration du terme d'erreur}
If we denote $f:=\left(-\frac{\partial^2}{\partial t^2}+(\vec{\Delta}+V)\right)de^{-t\sqrt{\Delta+V}} u$, then we have the following estimate on $f$:

$$||f(t,.)||_1+||f(t,.)||_p\leq \frac{C}{(1+t)^{n/p}}||u||_p.$$

\end{lem}
\textit{Proof of Lemma \ref{majoration du terme d'erreur}}\\\\
$$\begin{array}{lll}
\left(-\frac{\partial^2}{\partial t^2}+(\vec{\Delta}+V)\right)de^{-t\sqrt{\Delta+V}}u&=&-d(\Delta+V)e^{-t\sqrt{\Delta+V}}u+(\vec{\Delta}+V)de^{-t\sqrt{\Delta+V}}u \\\\
 
&=&-\big(e^{-t\sqrt{\Delta+V}}u\big)(dV).
\end{array}$$
We have:

$$||e^{-t\sqrt{\Delta+V}}||_{L^p\rightarrow L^\infty}\leq \frac{C}{t^{n/2p}},\,\forall t>0,$$
and

$$||e^{-t\sqrt{\Delta+V}}||_{L^p\rightarrow L^p}\leq 1,\,\forall t>0.$$
(this comes from the domination $e^{-t(\Delta+V)}\leq e^{-t\Delta}$). Since $V\in C_0^\infty$, we get the result.
\cqfd
Now we show:

\begin{lem}\label{Green estimates}

$||G(f)(t,.)||_2$ is bounded uniformly with respect to $t>0$, and 
$$\lim_{t\rightarrow0}||G(f)(t,.)||_2=0.$$

\end{lem}
\textit{Proof of Lemma \ref{Green estimates}:}\\\\
Denote $K_s(t,\sigma):=\frac{e^{-\frac{(\sigma-t)^2}{4s}}-e^{-\frac{(\sigma+t)^2}{4s}}}{\sqrt{4\pi s}}$.\\
$$\begin{array}{rcl}
G(f)(t,x)&=&\int G(\sigma,t,x,y)f(\sigma,y)d\sigma dy \\\\
&=&\int_M\int_0^\infty\int_0^\infty K_s(t,\sigma) \, \vec{p}^V_s(x,y)f(\sigma,y)\,dsd\sigma dy\\\\
&=&\int_0^\infty\int_0^\infty K_s(t,\sigma)\,\left(\int_M\vec{p}^V_s(x,y)f(\sigma,y)\,dy\right)\,dsd\sigma\\ \\
&=& \int_0^\infty\int_0^\infty K_s(t,\sigma)\,e^{-s\sqrt{\vec{\Delta}+V}}(x)\,dsd\sigma
\end{array} $$
Consequently,
$$||G(f)(t,.)||_2\leq\int_0^\infty\int_0^\infty K_s(t,\sigma)\, ||e^{-s\sqrt{\vec{\Delta}+V}}f(\sigma,.)||_2\, dsd\sigma.$$
But we have by Lemma \ref{ultracontractivite pour le Laplacien de Hodge},

$$\begin{array}{rcl}
||e^{-s\sqrt{\vec{\Delta}+V}}f(\sigma,.)||_2&\leq&\min\left(\frac{1}{s^{n/4}}||f(\sigma,.)||_1,||f(\sigma,.)||_2\right)\\\\
&\leq&C||u||_2\min\left(\frac{1}{s^{n/4}}\frac{1}{(1+\sigma)^{n/2}},\frac{1}{(1+\sigma)^{n/2}}\right)
\end{array}$$
Therefore, 

$$\begin{array}{rcl} ||G(f)(t,.)||_2&\leq& C||u||_2 \int_0^\infty\frac{1}{(1+\sigma)^{n/2}}\left(\int_0^1 \frac{e^{-\frac{(\sigma-t)^2}{4s}}-e^{-\frac{(\sigma+t)^2}{4s}}}{\sqrt{s}}\,ds\right.+\\\\
&&\hskip20mm\left.\int_1^\infty \frac{e^{-\frac{(\sigma-t)^2}{4s}}-e^{-\frac{(\sigma+t)^2}{4s}}}{s^{\frac{n}{4}+\frac{1}{2}}}\,ds\right)\,d\sigma
\end{array}$$
Since $n\geq3$, the three integrals $\int_0^\infty\frac{d\sigma}{(1+\sigma)^{n/2}}$, $\int_0^1\frac{ds}{\sqrt{s}}$ and $\int_1^\infty\frac{ds}{s^{\frac{n}{4}+\frac{1}{2}}}$ converge, and this yields immediately the fact that $ ||G(f)(t,.)||_2$ is bounded uniformly with respect to $t>0$. Furthermore, we can apply the Dominated Convergence Theorem to conclude that $\lim_{t\rightarrow0}||G(f)(t,.)||_2=0$.
\cqfd
\textit{End of the proof of Theorem \ref{transformee de Riesz avec potentiel}:}\\\\
Letting $(g(u))(x):=\int_0^\infty (G(f))(t,x)dt$, we have by integration of the formula in Proposition \ref{parametrix}:

$$(\vec{\Delta}+V)^{-1/2}du=d(\Delta+V)^{-1/2}u-cg.$$
By Lemma \ref{majoration du terme d'erreur} and Lemma \ref{ultracontractivite pour le Laplacien de Hodge}, we have as in \cite{Carron4}:

$$||g||_p\leq C||u||_p.$$
Applying Lemma \ref{(Delta+V)^-1/2 d borne}, we conclude that $d(\Delta+V)^{-1/2}$ is bounded on $L^p$, which finishes the proof of Theorem \ref{transformee de Riesz avec potentiel}.

\cqfd
\textit{Proof of Lemma \ref{ultracontractivite pour le Laplacien de Hodge}:}\\
Let us denote $L:=\vec{\Delta}+V$. If we can prove that $||e^{-tL}||_{\infty,\infty}\leq C$,  $||e^{-tL}||_{1,1}\leq C$ and $||e^{-tL}||_{1,\infty}\leq \frac{C}{t^{n/2}}$, then by standard interpolation arguments we are done. The fact that $||e^{-tL}||_{\infty,\infty}\leq C$ comes from the Gaussian estimate on $e^{-tL}$, which holds by Theorem \ref{estimee gaussienne pour le Laplacien de Hodge}, plus the fact that $\frac{1}{V(x,\sqrt{t})}\int_Me^{-c\frac{d^2(x,y)}{t}}dy$ is bounded uniformly in $x\in M$ and $t>0$. Then by duality $||e^{-tL}||_{1,1}\leq C$. Moreover, by Theorem \ref{estimee diagonale pour le Laplacien de Hodge} we also have the estimate:

$$||e^{-tL}||_{2,\infty}\leq \frac{C}{t^{n/4}},\,\forall t>0.$$
By duality, we deduce
$$||e^{-tL}||_{1,2}\leq \frac{C}{t^{n/4}},\,\forall t>0,$$
and by composition

$$||e^{-tL}||_{1,\infty}\leq||e^{-tL}||_{1,2}||e^{-tL}||_{2,\infty}\leq \frac{C^2}{t^{n/2}},\,\forall t>0.$$
The result for $e^{-t\sqrt{L}}$ follows by using the subordination identity:

$$e^{-t\sqrt{L}}=\frac{t}{2\sqrt{\pi}}\int_0^\infty e^{-\frac{t^2}{4s}}e^{-sL}\frac{ds}{s^{3/2}}.$$
\cqfd
\textit{Proof of Proposition \ref{parametrix}:}\\\\
Letting 
$$\varphi(t,.):=e^{-t\sqrt{\vec{\Delta}+V}}du-de^{-t\sqrt{\Delta+V}}u+G\left(\left(-\frac{\partial^2}{\partial t^2}+(\vec{\Delta}+V)\right)de^{-t\sqrt{\Delta+V}}u\right),$$
$\varphi(t,.)$ satisfies:

\begin{equation}\label{phi 1}
\left(-\frac{\partial^2}{\partial t^2}+(\vec{\Delta}+V)\right)\varphi=0,
\end{equation}
and also

\begin{lem}\label{phi 2}
The function $\varphi(t,.)$ converges to zero in $L^2$ when $t$ goes to $0$, that is

$$L^2-\lim_{t\rightarrow0}\varphi(t,.)=0.$$

\end{lem}
\pr Given the result of Lemma \ref{Green estimates} it is enough to prove that

\begin{equation}\label{egalite L2}
L^2-\lim _{t\rightarrow0}e^{-t\sqrt{\vec{\Delta}+V}}du=L^2-\lim _{t\rightarrow0}de^{-t\sqrt{\Delta+V}}u=du
\end{equation}
 In order to justify Equation \eqref{egalite L2}, we use first the Spectral Theorem ((c) in Theorem VIII.5 in \cite{Reed-Simon1}) which gives that $\sqrt{\Delta+V}e^{-t\sqrt{\Delta+V}}u$ converges in $L^2$ for $u\in C_0^\infty(M)$; since $V\geq0$, the Riesz transform with potential $d(\Delta+V)^{-1/2}$ is bounded on $L^2$, and from this we deduce that $de^{-t\sqrt{\Delta+V}}u$ converges in $L^2$; the limit is necessarily $du$.\cqfd
\textit{End of the proof of Proposition \ref{parametrix}:}\\\\
We claim that $\varphi(t,.)$ is bounded in $L^2$, uniformly with respect to $t>0$: to show this, it is enough to prove that $de^{-t\sqrt{\Delta+V}}u$ is uniformly bounded in $L^2$, when $u$ is a smooth, compactly supported fixed function. We write

$$de^{-t\sqrt{\Delta+V}}u=d\left(\Delta+V\right)^{-1/2}\sqrt{\Delta+V}e^{-t\sqrt{\Delta+V}}u,$$
and using the $L^2$ boundedness of $d\left(\Delta+V\right)^{-1/2}$ and the analyticity on $L^2$ of $e^{-t\sqrt{\Delta+V}}$, we get

$$||de^{-t\sqrt{\Delta+V}}u||_2\leq \frac{C}{t},\,\forall t>0.$$
Thus it remains to show that  $de^{-t\sqrt{\Delta+V}}u$ is bounded in $L^2$ when $t$ goes to $0$. But this follows from the fact that 

$$L^2-\lim_{t\rightarrow 0}de^{-t\sqrt{\Delta+V}}u=du,$$
which we have already proved. Consequently, $\varphi(t,.)$ is uniformly bounded in $L^2$. Using this, together with Equation \eqref{phi 1} and Lemma \ref{phi 2}, and applying the Spectral Theorem to $\vec{\Delta}+V$, we deduce that $\varphi\equiv 0$. This concludes the proof of Proposition \ref{parametrix}.\cqfd

\textit{This article is part of the PhD thesis of the author. He would like to express his gratitude for all the guidance and help provided by his advisor Prof. G. Carron. He would also to thank all the members of Laboratoire Jean Leray, Universit\'{e} de Nantes, for having provided him with a welcoming and stimulating working environment.}

\bibliographystyle{plain}

\bibliography{bibliographie}

\end{document}